\documentclass[12pt,oneside,reqno]{amsart}
\usepackage{amssymb}
\usepackage{color}
\usepackage{txfonts}
\usepackage{bbm}
\usepackage{cases}
\usepackage{amsmath}
\usepackage{graphicx}
\usepackage{mathrsfs}
\usepackage{stmaryrd}
\usepackage{amsfonts}
\usepackage{enumerate,amsmath,amssymb,amsthm}

\usepackage{color}
\usepackage[colorlinks=true]{hyperref}
\hypersetup{
    linkcolor=blue,          
    citecolor=red,        
    filecolor=blue,      
    urlcolor=cyan  
}

\numberwithin{equation}{section}

\allowdisplaybreaks

\newtheorem{theorem}{Theorem}[section]
\newtheorem{lemma}{Lemma}[section]

\newtheorem{proposition}{Proposition}[section]

\newtheorem{corollary}[theorem]{Corollary}
\newtheorem{remark}{Remark}[section]

\def\e{{\mathrm{e}}}
\def\eps{\varepsilon}

\def\a{\alpha}
\def\om{\omega}

\def\p{\partial}

\def\[{{\Big[}}
\def\]{{\Big]}}
\def\<{{\langle}}
\def\>{{\rangle}}
\def\({{\Big(}}
\def\){{\Big)}}

\def\bx{{\mathbf{x}}}

\def\dif{{\mathord{{\rm d}}}}

\def\min{{\mathord{{\rm min}}}}

\def\={&\!\!=\!\!&}
\def\bt{\begin{theorem}}
\def\et{\end{theorem}}
\def\bl{\begin{lemma}}
\def\el{\end{lemma}}
\def\br{\begin{remark}}
\def\er{\end{remark}}

\def\bd{\begin{definition}}
\def\ed{\end{definition}}
\def\bp{\begin{proposition}}
\def\ep{\end{proposition}}
\def\bc{\begin{corollary}}
\def\ec{\end{corollary}}
\def\bx{\begin{Examples}}
\def\ex{\end{Examples}}

\def\cB{{\mathcal B}}

\def\cE{{\mathcal E}}

\def\cI{{\mathcal I}}

\def\cN{{\mathcal N}}

\def\cP{{\mathcal P}}

\def\cT{{\mathcal T}}

\def\mC{{\mathbb C}}

\def\mE{{\mathbb E}}

\def\mH{{\mathbb H}}
\def\mI{{\mathbb I}}

\def\mL{{\mathbb L}}

\def\mN{{\mathbb N}}

\def\mP{{\mathbb P}}
\def\mQ{{\mathbb Q}}
\def\mR{{\mathbb R}}

\def\mT{{\mathbb T}}

\def\mZ{{\mathbb Z}}
\def\1{{\bf 1}}

\def\bP{{\mathbf P}}

\def\sA{{\mathscr A}}
\def\sB{{\mathscr B}}

\def\sD{{\mathscr D}}

\def\sF{{\mathscr F}}
\def\sG{{\mathscr G}}

\def\sI{{\mathscr I}}

\def\sL{{\mathscr L}}
\def\sM{{\mathscr M}}

\def\sP{{\mathscr P}}

\def\sV{{\mathscr V}}

\def\geq{\geqslant}
\def\leq{\leqslant}

\def\div{\mathord{{\rm div}}}

\def\nor{|\mspace{-1mu}|\mspace{-1mu}|}

\def\u{{\mathbf u}}

\def\bP{{\mathbf P}}

\def\bE{{\mathbf E}}

\newtheorem{definition}{Definition}

\def\mN{{\mathbb N}}
\def\mR{{\mathbb R}}

\begin{document}

\title{Stochastic Lagrangian path for Leray solutions of 3D Navier-Stokes equations}
\author{Xicheng Zhang and Guohuan Zhao}

\address{Xicheng Zhang:
School of Mathematics and Statistics, Wuhan University,
Wuhan, Hubei 430072, P.R.China\\
Email: XichengZhang@gmail.com
 }

\address{Guohuan Zhao:
Fakult\"at f\"ur Mathematik, Universit\"at Bielefeld,
33615, Bielefeld, Germany\\
Email: zhaoguohuan@gmail.com
 }

\thanks{
Research of X. Zhang is partially supported by NNSFC grant of China (No. 11731009). Research of G. Zhao is supported by the German Research Foundation (DFG) through the Collaborative Research Centre(CRC) 1283 ``Taming uncertainty and profiting from randomness and low regularity in analysis, stochastics and their applications".
}

\begin{abstract}
In this paper we show the existence of stochastic Lagrangian particle trajectory for Leray's solution of 3D Navier-Stokes equations.
More precisely, for any Leray's solution ${\mathbf u}$ of 3D-NSE and each $(s,x)\in\mathbb{R}_+\times\mathbb{R}^3$, we show the existence of weak solutions to the following SDE, which has a density $\rho_{s,x}(t,y)$ belonging to $\mathbb{H}^{1,p}_q$ provided $p,q\in[1,2)$ with $\frac{3}{p}+\frac{2}{q}>4$:
$$
\mathrm{d} X_{s,t}={\mathbf u} (s,X_{s,t})\mathrm{d} t+\sqrt{2\nu}\mathrm{d} W_t,\ \ X_{s,s}=x,\ \ t\geq s,
$$
where $W$ is a three dimensional standard Brownian motion, $\nu>0$ is the viscosity constant. Moreover,
we also show that for Lebesgue almost all $(s,x)$, the solution $X^n_{s,\cdot}(x)$ of the above SDE associated with the 
mollifying velocity field ${\mathbf u}_n$ weakly converges to $X_{s,\cdot}(x)$ so that $X$ is a Markov process in almost sure sense. 


\end{abstract}

\maketitle \rm


\section{Introduction}

Throughout the paper we assume $d\geq 2$. Consider the following Navier-Stokes equation:
$$
\p_t\u=\nu\Delta\u+\u\cdot\nabla \u+\nabla p,\ \ \div\,\u\equiv0,\ \u_0=\varphi,
$$
where $\u=(u_1,\cdots, u_d)$ is the velocity field of the fluid, $\nu>0$ is the viscosity constant, and $p$ stands for the pressure.
It is well known that for any divergence free vector field $\varphi\in L^2(\mR^d)$, there exists a divergence free Leray weak solution to 3D-NSEs in the class 
\begin{align}\label{FD1}
\|\u\|_{L^\infty([0,T]; L^2(\mR^d))}+\|\nabla\u\|_{L^2([0,T]; L^2(\mR^d))}<\infty,\ \ \forall T>0.
\end{align}
In a recent remarkable paper,  Buckmaster and Vicol \cite{Bu-Vi} showed that
there are infinitely many weak solutions $\u\in C(\mR_+; L^2(\mT^3))$ for 3D-NSEs on the torus.
However, it is still not known whether the above Leray solution is unique and smooth, which are in fact famous 
open problems for a long time.  

\medskip

In this work we are interesting in the following problem:
For any Leray solution $\u$, is it possible to construct the stochastic Lagrangian particle trajectory $X_t=X_t(x)$ 
associated with the velocity field $\u$? More precisely, for each starting point $x$, is there a unique solution to the following SDE?
\begin{align}\label{SDE0}
\dif X_t=\u(t,X_t)\dif t+\sqrt{2\nu}\dif W_t,\ X_0=x,
\end{align}
where $W$ is a $d$-dimensional standard Brownian motion on some probability space $(\Omega,\sF,\bP)$.
If $\u$ is smooth in $x$,  then  by Constantin and Iyer's representation \cite{Co-Iy} (see also \cite{Zh2, Zh3}), $\u$ can be reconstructed from $X_t(x)$ as follows:
$$
\u(t,x)=\cP\mE(\nabla^{\rm t} X^{-1}_t(x)\cdot\varphi(X^{-1}_t(x))),
$$
where $\cP$ is the Leray projection and $X^{-1}_t(x)$ is the inverse of stochastic flow $x\mapsto X_t(x)$,
and $\nabla^{\rm t}$ stands for the transpose of Jacobian matrix.
By Krylov and R\"ockner's result \cite{Kr-Ro},  under the following assumption
\begin{align}\label{Red}
\u\in  \cap_{T>0}L^q([0,T]; L^p(\mR^d)),\ \ p,q\geq 2,\ \ \tfrac{d}{p}+\tfrac{2}{q}<1,
\end{align}
there is a unique strong solution to SDE \eqref{SDE0} for any starting point $x\in\mR^d$. Moreover, the unique solution $X_t(x)$ 
is weakly differentiable in $x$ and satisfies (see \cite{Fe-Fl, Zh0, Zh1}):
$$
\sup_{x\in\mR^d}\bE\left(\sup_{t\in[0,T]}|\nabla X_t(x)|^p\right)<\infty,\ \ \forall p\geq 1,\ T>0.
$$
Unfortunately, Leray's solution does not satisfy \eqref{Red}. Indeed, by \eqref{FD1} and Sobolev's embedding (see Lemma \ref{Le11} below),
\begin{align}\label{LK0}
\u\in \cap_{T>0}L^q([0,T]; L^p(\mR^d)),\ \ p,q\geq 2,\ \ 
\tfrac{d}{p}+\tfrac{2}{q}>\tfrac{d}{2}.
\end{align}
Notice that the {\it deterministic} Lagrangian particle trajectories associated with $\u$  
have been studied very well (for example, see \cite[Chapter 17]{Ro-Ro-Sa} and \cite{Ch-Le}),
which depends on further regularity on Leray's solution. Here we want to solve SDE \eqref{SDE0} only basing on \eqref{LK0} for $d=3$.

\medskip

For given $(s,x)\in\mR_+\times\mR^d$, we consider the following SDE in $\mR^d$ starting from $x$ at time $s$:
\begin{align}\label{SDE}
\dif X_{s,t}=b(t,X_{s,t})\dif t+\sqrt{2}\dif W_t,\ t>s, \ X_{s,s}=x,
\end{align}
where $b(t,x): \mR_+\times\mR^d\to\mR^d$ is a measurable vector field.
The generator associated with the above SDE is given by
$$
\sL^b_t:=\Delta+b(t,\cdot)\cdot\nabla.
$$
 In this paper, we focus on the weak solution of SDE \eqref{SDE} with {\it lower regularity} $b$, that is,
 $$
 b\in \cap_{T>0}L^q([0,T];L^p(\mR^d))=:L^q_{loc}(L^p),\ \ p,q\geq 2,\ \ \tfrac{d}{p}+\tfrac{2}{q}<2.
 $$
Roughly speaking, a weak solution of SDE \eqref{SDE} is a {\em semimartingale} $(X_{s,t})_{t\geq s}$ so that 
$$
\int_s^t |b(r,X_{s,r})| \dif r<\infty,\ \ \forall t\geq s,\quad a.s.,
$$
and
\begin{align}\label{SDE2}
X_{s,t}=x+ \int_s^t b(r,X_{s,r})\dif r+\sqrt{2}(W_t-W_s), \ \ \forall t\geq s,\quad a.s. 
\end{align}

\medskip

When $b\in L^q_{loc}(L^p)$ for some $p,q\in[2,\infty)$ with $\frac{d}{p}+\frac{2}{q}<1$, as mentioned above,
by Girsanov's transformation and $L^p$-theory of second order parabolic equations, 
Krylov and R\"ockner \cite{Kr-Ro} showed that there is a unique strong solution to SDE \eqref{SDE}, 
which extended the main results in \cite{Ve} and \cite{Zv}. 
The strong well-posedness of SDE \eqref{SDE} driven by multiplicative Brownian noise  
was studied in \cite{Zh0, Zh1} by Zvonkin's transformation introduced in \cite{Zv}. 
Moreover, the flow property and weak differentiability of $X_{s,t}(x)$ in $x$ are also obtained therein.
When $b\in H^{-\alpha,p}$ with $\alpha\in(0,\frac{1}{2})$ and $p\in(\frac{d}{1-\alpha},\frac{d}{\alpha})$ is time-independent,
Flandoli, Issoglio and Russo \cite{Fl-Is-Ru} showed the existence and uniqueness of
``virtual''  solutions (a class of special weak solutions) to SDE \eqref{SDE}. Later, the well-posedness of martingale solutions and weak solutions
(which may not be a semimartingale but a Dirichlet process) was established in \cite{Zh-Zh} 
for $b\in H^{-\a,p}$ with $\a\in (0,\frac{1}{2}]$ and $p\in (\frac{d}{1-\a},\infty)$. 
We also mention that 
Bass and Chen  in \cite{Ba-Ch2} studied the weak well-posedness of SDE \eqref{SDE} in the class of semimartingales 
when $b$ belongs to some generalized Kato's class ${\bf K}_{d-1}$(see also \cite{Zh-Zh1}), in particular, the space $L^p$ with $p>d$ is included in that class. 

\medskip

It should be emphasized that even in the weak sense, all the works mentioned above does not cover the borderline 
case $b\in L^q_{loc}(L^p)$ with $\frac{d}{p}+\frac{2}{q}=1$,  
not to mention the supercritical case 
$\frac{d}{p}+\frac{2}{q}>1$.  Let us explain the difficulty firstly. 
In order to get the weak existence of SDE \eqref{SDE} with singular drifts, a straightforward way is to use Girsanov's transform as in \cite{Kr-Ro}. 
However, this approach does not work in the case when $p\leq d$.  Let us make a detailed analysis for this point.
Let $\mC$ be the space of all continuous functions from $\mR_+$ to $\mR^d$, which is endowed with the usual Borel $\sigma$-field $\cB(\mC)$.
All the probability measures over $(\mC,\cB(\mC))$ is denoted by $\sP(\mC)$. 
Let $\omega_t$ be the canonical process over $\mC$. For $t\geq 0$, let $\cB_t:=\cB_t(\mC)$ be the natural filtration generated by $\{\omega_s: s\leq t\}$.
Let $\mP\in\sP(\mC)$ be the classical Wiener measure so that $t\mapsto\omega_t$ is a $d$-dimensional standard Brownian motion.
For $b\in L^p(\mR^d)$ with $p\leq d$, one can check that the Novikov condition 
\begin{align}\label{Con1}
\mE \exp\left(\frac{1}{2}\int_0^T |b|^2(\omega_t) \dif t\right)<\infty
\end{align}
for the exponential supermartingale 
$$
\cE^b_t=\exp\left(\int_0^t b(\omega_s)\dif \omega_s-\frac{1}{2}\int_0^t b(\omega_s) \dif s\right) 
$$ 
may not hold. Notice that condition \eqref{Con1} is somehow  equivalent 
to say that $b$ belongs to some Kato's class (see \cite{Ai-Si}).
In fact, without other conditions, if $b$ only belongs to $L^{d-\eps}_{loc}\setminus L^d_{loc}$, then the weak existence may be failed.
For example, consider the following SDE:
\begin{align}\label{SDE09}
X_t =-c\int^t_0 X_s|X_s|^{-2}\dif s + W_t,\ \ c\in\mR.
\end{align}
If $c\geq d$, Kinzebulatov and Semenov \cite[page 3]{Ki-Se} explained why the above SDE does not allow a solution (see also \cite{Be-Fl-Gu-Ma}).
Meanwhile, for $c<c_d$, where $c_d\in(0,d)$ is some constant only depending on $d$, they proved that there exists a weak solution to the above SDE
by utilizing the analytic construction of the semigroup $\e^{-t(\Delta+b\cdot\nabla)}$.
By direct calculations, for $b(x):=-c x|x|^{-2}$ and $d\geq 3$, we have
$$
\div b(x)=-c(d-2)|x|^{-2}\notin L_{loc}^{d/2}.
$$
Intuitively, if $c\geq d$, then the centripetal force is so strong such that the particle can not escape from the origin immediately 
so that even though a random perturbation is added, there is no solution for SDE \eqref{SDE09}. 
However, our result below shows that if $b\in L^{d/2+\eps}(\mR^d)$ for some $\eps>0$, then
equation \eqref{SDE} has at least one semimartingale solution,  provided that the 
negative part of $\div b$ satisfies some integrable conditions. 
We emphasize that Kinzebulatov and Semenov's result in \cite{Ki-Se} can not be applied to the case $b \in L^{d-\eps}_{loc}\setminus L^d_{loc}$.
We believe that the divergence condition is necessary for this case. Moreover, 
the singular time-dependent drift $b$ is not treated in \cite{Ki-Se}. If it is not possible, it seems hard to directly construct the two-parameter semigroups
associated with time-dependent drifts by analytic method.

\medskip

Before stating our results, we introduce the following notion of martingale solutions.
\bd
For given $(s,x)\in\mR_+\times\mR^d$, we call a probability measure $\mP_{s,x}\in\sP(\mC)$ a martingale solution of SDE \eqref{SDE} with starting point $(s,x)$ if
\begin{enumerate}[(i)]
\item $\mP_{s,x}(\omega_t=x,t\leq s)=1$, and for each $t>s$,
$$
\mE^{\mP_{s,x}}\left(\int^t_s |b(r,\omega_r)|\dif r\right)<\infty.
$$
\item For all $f\in C^2_c(\mR^d)$, $M^f_t$ is a $\cB_t$-martingale under $\mP_{s,x}$, where
$$
M^f_t(\omega):=f(\omega_t)-f(x)-\int^t_s\sL^b_rf(\omega_r)\dif r,\ \ t\geq s.
$$
\end{enumerate}
All the martingale solution $\mP_{s,x}$ with starting point $(s,x)$ and drift $b$ is denoted by $\sM^b_{s,x}$.
\ed
\br
Let $\mP_{s,x}\in\sM^b_{s,x}$. By L\'evy's characterization for Brownian motion, one sees that
$$
W_t:=\frac{\sqrt{2}}{2}\left(\omega_t-\omega_s-\int^t_sb(r,\omega_r)\dif r\right), \ \ t\geq s,
$$
is a $d$-dimensional standard Browian motion under $\mP_{s,x}$ (see \cite[Theorem 4.2.1]{St-Va}), so that 
$$
\omega_t=x+\int^t_sb(r,\omega_r)\dif r+\sqrt{2}W_t,\ \ t\geq s.
$$
In other words, $(\mC,\sB(\mC),\mP_{s,x}, \omega_t, W_t)$ is a weak solution of SDE \eqref{SDE}.
\er

Our main result is
\bt\label{TH11}
 Suppose that for some $p_i,q_i\in[2,\infty)$ with $\tfrac{d}{p_i}+\tfrac{2}{q_i}<2$, $i=1,2$,
\begin{align}\label{Con2}
\nor b\nor_{0,p_1;q_1}+\nor(\div\, b)^-\nor_{0,p_2;q_2}<\infty,
\end{align}
where $\nor\cdot\nor_{\alpha, p;q}$ is defined by \eqref{FQ2} below.
For each $(s,x)\in\mR_+\times\mR^d$, there exists at least one martingale solution $\mP_{s,x}\in\sM^b_{s,x}$, 
which satisfies the following Krylov's type estimate:
for any $\alpha\in[0,1]$ and $p,q\in(1,\infty)$ with $\frac{d}{p}+\frac{2}{q}<2-\alpha$, there exist $\theta=\theta(\alpha,p,q)>0$
and a constant $C>0$
such that for all $s\leq t_0<t_1<\infty$ with $t_1-t_0\leq 1$ and $f\in C^\infty_c(\mR^{d+1})$,
\begin{align}\label{EM9}
\mE^{\mP_{s,x}}\left(\int^{t_1}_{t_0} f(t,\omega_t)\dif t\Bigg|\cB_{t_0}\right)\leq C(t_1-t_0)^\theta\nor f\nor_{-\alpha,p;q}.
\end{align}
Moreover, we have the following conclusions:
\begin{enumerate}[(i)]
\item (Weak uniqueness)  For any mollifying approximation  $b_n$ of $b$,
there is a Lebesgue-null set $\cN\subset\mR_+\times\mR^d$ such that for all $(s,x)\in\cN^c$,
\begin{align}\label{JH1}
\mP^n_{s,x}\mbox{ weakly converges to $\mP_{s,x}\in\sM^b_{s,x}$}, 
\mbox{ where $\mP^n_{s,x}\in\sM^{b_n}_{s,x}$}.
\end{align}
\item (Almost surely Markov property) For each $(s,x)\in\cN^c$, there is a Lebesgue null set $I_{s,x}\subset[s,\infty)$ such that 
for all $t_0\in (s,\infty)\setminus I_{s,x}$, any $t_1>t_0$ and $f\in C_c(\mR^d)$,
\begin{align}\label{JH2}
\mE^{\mP_{s,x}}(f(\omega_{t_1})|\cB_{t_0})=\mE^{\mP_{t_0,\omega_{t_0}}}(f(\omega_{t_1})),\ \ \mP_{s,x}-a.s.
\end{align}
\item ($L^p$-semigroup)
Let $\cT_{s,t}f(x):=\mE^{\mP_{s,x}} f(\omega_t)$.
For any $p\geq 1$ and $T>0$, there is a constant $C>0$ such that for Lebesgue almost all $0\leq s<t\leq T$ and $f\in L^p(\mR^d)$,
\begin{align}\label{Co1}
\|\cT_{s,t}f\|_p\leq C\|f\|_p.
\end{align}
\end{enumerate}
\et

\br
By discretization stopping time approximation, Krylov estimate \eqref{EM9} is equivalent to say that 
for any $\delta\in(0,1)$ and stopping time $\tau\in[s,\infty)$,
\begin{equation}\label{Eq-Kry2}
\mE^{\mP_{s,x}} \left(\int_{\tau}^{\tau+\delta} f(t,\om_t)\dif t  \Big|\cB_{\tau} \right) \leq C\delta^\theta\nor f\nor_{\alpha, p;q},
\end{equation}
where $\cB_\tau:=\sigma\big\{\omega_{t\wedge\tau}, t\geq 0\big\}$ is the stopping time $\sigma$-field. In fact, let $\tau_n$ 
be a sequence of decreasing stopping times taking values in $\mT:=\{k\cdot 2^{-n}: k, n\in \mN\}$ and converging to $\tau$.
For any $f\in C^\infty_c(\mR^{d+1})$ and $\delta\in(0,1)$, by the dominated convergence theorem and martingale convergence theorem, we have
\begin{align*}
&\mE^{\mP_{s,x}} \left(\int_{\tau}^{\tau+\delta} f(t,\om_t)\dif t  \Big|\cB_{\tau} \right)
=\lim_{n\to\infty}\mE^{\mP_{s,x}} \left(\int_{\tau_n}^{\tau_n+\delta} f(t,\om_t)\dif t  \Big|\cB_{\tau_n} \right)\\
&\qquad=\lim_{n\to\infty}\mE^{\mP_{s,x}}\left(\sum_{a\in\mT}\1_{\{\tau_n=a\}}
\int_{a}^{a+\delta} f(t,\om_t)\dif t  \Big|\cB_{\tau_n} \right)\\
&\qquad=\lim_{n\to\infty}\sum_{a\in\mT}\1_{\{\tau_n=a\}}
\mE^{\mP_{s,x}} \left(\int_a^{a+\delta} f(t,\om_t)\dif t  \Big|\cB_{a} \right)\\
&\qquad\leq C\delta^\theta\nor f\nor_{-\alpha, p;q}\lim_{n\to\infty}\sum_{a\in\mT}\1_{\{\tau_n=a\}}=C\delta^\theta\nor f\nor_{-\alpha, p;q}.
\end{align*}
Moreover, let $\mu_{s,x}(t,\dif y):=\mP_{s,x}\circ \omega_t^{-1}$. 
For any $\alpha\in[0,1]$ and $p,q\in(1,\infty)$ with $\frac{d}{p}+\frac{2}{q}<2-\alpha$, by \eqref{EM9}, for any $T>0$ there is a constant $C>0$ such that
for all $f\in C^\infty_c([0,T]\times\mR^d)$,
$$
\left|\int^T_0\!\!\int_{\mR^{d}}f(t,y)\mu_{s,x}(t,\dif y)\dif t\right|\leq C\|f\|_{-\alpha,p;q},
$$
which in turn implies that $\mu_{s,x}(t,\dif y)=\rho_{s,x}(t,y)\dif y$ with $\rho_{s,x}\in\mH^{\alpha, p/(p-1)}_{q/(q-1)}$.
\er

\br
If $(\div b)^-\equiv 0$, then $\|\cT_{s,t}f\|_1\leq\|f\|_1$ in \eqref{Co1}.
If $\div b\equiv0$, then for any nonnegative $f\in L^1(\mR^d)$,
 $\|\cT_{s,t}f\|_1=\|f\|_1$.
By \eqref{LK0}, we can apply the above theorem to the Leray solution of 3D-NSEs. 
\er

\br
Let $d\geq 3$ and $\alpha<3$. Define
$$
b(x):=\sum_{z\in\mZ^d}\gamma_z \frac{x-z}{|x-z|^{\alpha}}\phi(|x-z|),
$$
where for some $M>0$, $\gamma_z\in(0,M)$ is a constant and $\phi\in C^\infty_c(\mR_+;[0,1])$ with $\phi(r)=1$ for $r\in[0,1]$ and
$\phi(r)=0$ for $r>2$. It is easy to see that \eqref{Con2} holds.
\er

\br
It should be compared with the results in \cite{Zh23, Zh31}. Therein, under the assumptions 
\begin{align}\label{Con9}
\nabla b\in \mL^1_{loc},\quad (\div b)^-, b/(1+|x|)\in\mL^\infty,
\end{align}
the existence and uniqueness of almost everywhere stochastic flows are obtained in the framework of DiPerna-Lions' theory.
By the estimate \eqref{Co1}, we can weaken the assumption on the boundedness of $(\div b)^-$ in \cite{Zh31} when the noise is nondegenerate. 
On the other hand, in \cite{Zh23, Zh31},
under \eqref{Con9},
the existence of a solution is only shown for Lebesgue almost all $x\in\mR^d$, while, under \eqref{Con2} we can show the existence of a solution for all starting point $x\in\mR^d$.
\er

To prove Theorem \ref{TH11}, the key point for us is to establish the maximum principle for the following parabolic equation under \eqref{Con2}:
\begin{align}\label{FD2}
\p_t u=\Delta u+b\cdot \nabla u+f,\  u(0)=0.
\end{align}
More precisely, for any $\alpha\in[0,1]$ and $q,p\in(1,\infty)$ with $\frac{d}{p}+\frac{2}{q}<2-\alpha$,
\begin{align}\label{Ma}
\|u\|_{L^\infty([0,T]\times\mR^d)}\leq C\nor f\nor_{-\alpha,p;q}.
\end{align}
When $f\equiv 0$,  under \eqref{Con2} the local maximum principle is proved by Nazarov and Ural'tseva in \cite{Na-Ur} by using Moser's iteration.
We also refer to \cite{Ha-Li} for the study of elliptic equations with drift $b\equiv 0$ and $f\in L^p(\mR^d)$ for $p>\frac{d}{2}$.
Here an open question is that whether we can show \eqref{JH1}-\eqref{Co1} for all $(s,x)\in\mR_+\times\mR^d$, 
which is closely related to find a continuous solution for PDE \eqref{FD2} under \eqref{Con2}. 

\medskip

This paper is organized as follows: 
In Section 2, we establish the key maximum estimate \eqref{Ma} by De Giorgi's method. In fact, we shall show a more general result by 
allowing $b$ and $f$ being in negative Sobolev spaces, which are not treated in \cite{Ha-Li, Na-Ur}. 
In Section 3, we prove our main result Theorem \ref{TH11}. In Appendix, some properties of certain
local Sobolev spaces are given.
Throughout this paper we shall use the following conventions:
\begin{itemize}
\item We use $A\lesssim B$  to denote $A\leq C B$ for some unimportant constant $C>0$.  
\item For any $\eps\in(0,1)$, we use $A\lesssim\eps B+D$ to denote
$A\leq \eps B+C_\eps D$ for some constant $C_\eps>0$.
\item $\mN_0:=\mN\cup\{0\}$, $\mR_+:=[0,\infty)$, $a\vee b:=\max(a,b)$, $a\wedge b:=\min(a,b)$, $a^+:=a\vee 0$.

\end{itemize}

\section{Maximum principle for parabolic equations by De Giorgi's method}

We first introduce some spaces, functions and notations for later use:
\begin{itemize}
\item Let $C^\infty_c(\mR^{d+1})$ be the space of all smooth functions with compact supports and $\sD'$ the dual space of 
$C^\infty_c(\mR^{d+1})$, which is also called distribution space. 
The duality between $\sD'$ and $C^\infty_c(\mR^{d+1})$ is denoted by $\<\!\<\cdot,\cdot\>\!\>$.
In particular, if $f(t,x)$ and $g(t,x)$ are two real functions in $\mR\times\mR^{d}$, then
$$
\<\!\<f,g\>\!\>=\int_\mR\<f(t),g(t)\>\dif t\ \ \mbox{with\ \ } \<f(t),g(t)\>:=\int_{\mR^d}f(t,x)g(t,x)\dif x.
$$
\item For two distributions $f,g\in\sD'$, one says that $f\leq g$ if for any nonnegative $\varphi\in C^\infty_c(\mR^{d+1})$,
$$
\<\!\<f,\varphi\>\!\>\leq \<\!\<g,\varphi\>\!\>.
$$

\item For $\alpha\in\mR$ and $p\in[1,\infty]$, let $H^{\alpha,p}$ be the usual Bessel potential space with norm:
$$
\|f\|_{\alpha,p}:=\|(\mI-\Delta)^{\alpha/2}f\|_p=\left(\int_{\mR^d}|(\mI-\Delta)^{\alpha/2}f(x)|^p\dif x\right)^{1/p}.
$$
If $f\in H^{\a, p}$, $g\in H^{-\a, p'}$ with $\frac{1}{p}+\frac{1}{p'}=1$, we denote 
$$
\<f, g\>:= \int_{\mR^d} (\mI-\Delta)^{-\a/2}f (x)\cdot (\mI-\Delta)^{\a/2} g(x) \dif x.
$$
\item For $\alpha\in\mR$ and $p,q\in[1,\infty]$, let $\mH^{\alpha,p}_q:=L^q(\mR; H^{\alpha,p})$ be the space of spatial-time functions with norm
$$
\|f\|_{\alpha,p;q}:=\left(\int_{\mR}\|f(t,\cdot)\|_{\alpha,p}^q\dif t\right)^{1/q}.
$$
If $f\in \mH^{\a, p}_{q}$, $g\in \mH^{-\a, p'}_{q'}$ with $\frac{1}{p}+\frac{1}{p'}=1$, $\frac{1}{q}+\frac{1}{q'}=1$, we also denote 
$$
\<\!\<f, g\>\!\>:= \int_{\mR} \<f(t), g(t)\> \dif t=\int_{\mR} \<f, g\>. 
$$
\item Let $\mH^{\alpha,p}_{q,loc}$ be the space of all functions $f:\mR^{d+1}\to\mR$ with
$$
f\eta\in\mH^{\alpha,p}_q,\ \forall \eta\in C^\infty_c(\mR^{d+1}).
$$
\item For $r>0$, we define 
$$
B_r:=\{x\in\mR^d: |x|<r\},\ Q_r:=(-r^2,r^2)\times B_r.
$$ 
\item Fix $\chi\in C^\infty(\mR^{d+1};[0,1])$ with $\chi|_{Q_1}=1$ and $\chi|_{Q^c_2}=0$. For $r>0$ and $(s,z)\in\mR^{d+1}$, define
\begin{align}\label{CH}
\chi_r(t,x):=\chi(r^{-2}t,r^{-1}x),\ \chi^{s,z}_r(t,x):=\chi_r(t-s,x-z),\ \ (t,x)\in\mR^{d+1}.
\end{align}
\item Fix $r>0$. Let $\widetilde\mH^{\alpha,p}_q$ be the Banach space of all functions $f\in \mH^{\alpha,p}_{q,loc}$ with
\begin{align}\label{FQ2}
\nor f\nor_{\alpha,p;q}:=\sup_{s,z}\|f\chi^{s,z}_r\|_{\alpha,p;q}<\infty.
\end{align}
\item We shall simply write
$$
\mL^p_q:=\mH^{0,p}_q,\ \ \|f\|_{p;q}:=\|f\|_{0,p;q},\ \ \widetilde\mL^p_q:=\widetilde\mH^{0,p}_q,\ \ \nor f\nor_{p;q}:=\nor f\nor_{0,p;q}.
$$
\item Let $\sV:=\mL^2_\infty\cap\mH^{1,2}_2$, $\sV_{loc}:=\mL^2_{\infty,loc}\cap\mH^{1,2}_{2,loc}$
and $\widetilde\sV:=\widetilde\mL^2_\infty\cap\widetilde\mH^{1,2}_{2}$. Define
\begin{align}\label{NK1}
\|f\|_{\sV}:=\|f\|_{2;\infty}+\|\nabla_xf\|_{2;2},\ \ \nor f\nor_{\sV}
:=\nor f\nor_{2;\infty}+\nor\nabla_xf\nor_{2;2}.
\end{align}
\item Let $\rho\in C^\infty_c(B_1)$ with $\int\rho=1$. For $\eps\in(0,1)$ and $n\in\mN$, we shall use the mollifiers:
$$
\rho_\eps(x):=\eps^{-d}\rho(\eps^{-1} x),\ \ \rho_n(x):=\rho_{1/n}(x).
$$
\item For $d\geq 2$, define
$$
\sI_d:=\left\{(\alpha,p,q)\in[0,1]\times(1,\infty)\times(1,\infty): \frac{d}{p}+\frac{2}{q}<2-\alpha\right\}.
$$
\item For given $(\alpha,p,q)\in\sI_d$, we define $r,s\in[2,\infty]$ by relation
\begin{align}\label{NH19}
\frac{1}{(2-\alpha)p}+\frac{1}{r}=\frac{1}{(2-\alpha)q}+\frac{1}{s}=\frac{1}{2}.
\end{align}
Notice that $\frac{d}{p}+\frac{2}{q}<2-\alpha$ implies
$$
\frac{d}{r}+\frac{2}{s}>\frac{d}{2}.
$$
\item The following Gagliardo-Nirenberge's interpolation inequality will be used frequently:
\begin{align}\label{EQ1}
\|u\|_{\alpha,r}\leq C\|\nabla u\|^\theta_p\|u\|^{1-\theta}_q,
\end{align}
where $\alpha\in[0,1]$, $\theta\in[\alpha,1]$ and $p,q,r\in[1,\infty]$ satisfy
$$
\frac{1}{r}=\frac{\alpha}{d}+\theta\left(\frac{1}{p}-\frac{1}{d}\right)+\frac{1-\theta}{q},\ \ 1-\frac{d}{p}\notin\mN_0.
$$
\end{itemize}
\subsection{Localization estimates}
In this subsection we prove an important localization lemma for later use, which is a consequence of Gagliado-Nirenberge's interpolation inequality and H\"older's inequality.
First of all, we have the following interpolation estimates.
\bl\label{Le11}
Let $r,s\geq2$ with $\frac{d}{r}+\frac{2}{s}>\frac{d}{2}$. 
For any $\eps\in(0,1)$, there is a constant $C_\eps>0$ such that
$$
\|f\|_{r;s}\leq \eps\|\nabla f\|_{2;2}+C_\eps\|f\|_{2;2(1-\theta)s/(2-s\theta)},\quad \forall f\in \mH^{1,2}_2\cap\mL^2_{2(1-\theta)s/(2-s\theta)},
$$
where $\theta:=\frac{d}{2}-\frac{d}{r}\in [0,1)$.
In particular, if {\rm supp}$f\subset Q_2$, then
for some $C=C(d,r,s)>0$,
\begin{align}\label{EB}
\|f\|_{r;s}\leq C\|f\|_\sV,\ f\in\sV.
\end{align}
\el
\begin{proof}
For $r\in [2,\infty)$ if $d=2$ or $r\in[2,2d/(d-2)]$ if $d\geq 3$,
by \eqref{EQ1} we have
$$
\|f\|_r\leq C\|\nabla f\|^{\theta}_2\|f\|^{1-\theta}_2.
$$
Since $s\theta<2$, by H\"older's inequality we further have
$$
\|f\|_{r;s}\leq C\|\nabla f\|^{\theta}_{2;2}\|f\|^{1-\theta}_{2;2(1-\theta)s/(2-s\theta)},
$$
which gives the desired embedding by Young's inequality.
\end{proof}
\bl\label{Le12}
Let $Q=I\times D\subset \mR\times\mR^{d}$ 
be a bounded domain. For any $p,q,r,s\in[1,\infty]$, there is a constant $C>0$ only depending on $Q, p,q,r,s$ such that for any $A\subset Q$,
$$
\|\1_A\|_{p;q}\leq C\|\1_A\|^{(r/p)\wedge(s/q)}_{r;s}.
$$
\el
\begin{proof}
Define
$$
A_t:=\int_D \1_A(t,x)\dif x.
$$
If $r/p\leq s/q$, then by H\"older's inequality,
$$
\|\1_A\|_{p;q}=\left(\int_I A_t^{q/p}\dif t\right)^{1/q}\leq C \left(\int_I A_t^{s/r}\dif t\right)^{r/(sp)}=C\|\1_A\|^{r/p}_{r;s}.
$$
If $r/p>s/q$, then by H\"older's inequality,
$$
A_t^{q/p}\leq C A_t^{s/r}\Rightarrow \|\1_A\|_{p;q}=\left(\int_I A_t^{q/p}\dif t\right)^{1/q}\leq 
 C\left(\int_I A_t^{s/r}\dif t\right)^{1/q}=C\|\1_A\|^{s/q}_{r;s}.
$$
The proof is complete.
\end{proof}
The following lemma is the key localization result. 
\bl
Let $\eta:\mR^{d+1}\to[0,1]$ be a smooth function with compact support contained in $Q_2$.
Let $(\alpha,p,q)\in\sI_d$ and $r,s\in[2,\infty]$ be defined by \eqref{NH19}. Let $\chi_2$ be the cutoff function defined by \eqref{CH}.
For any $\eps\in(0,1)$, there is a constant $C_\eps=C_\eps(d,\alpha,p,q)>0$ 
such that for any $c, b,f\in\mH^{-\alpha,p}_{q,loc}$ and $w\in\sV_{loc}$,
\begin{align}\label{NJ52}
|\<\!\<c, \eta^2w^2\>\!\>|
\leq\eps\|\eta w\|^2_{\sV}+C_\eps\|c\chi_2\|^2_{-\alpha,p;q} \|\eta w\|_{r;s}^2,
\end{align}
\begin{align}\label{NJ51}
|\<\!\< b, \nabla\eta^2w^2\>\!\>|
\leq\eps\|\eta w\|^2_{\sV}+C_\eps\Big(1+\|b\chi_2\|^2_{-\alpha,p;q}\Big)
\Big(1+\|\nabla^2\eta\|_\infty+\|\nabla\eta\|_\infty^2\Big)\|w\1_{\eta\not=0}\|_{r;s}^2,
\end{align}
\begin{align}\label{NJ55}
|\<\!\<f, \eta^2w\>\!\>|
\leq\eps\|\eta w\|^2_{\sV}+C_\eps\|f\chi_2\|^2_{-\alpha,p;q}\Big(1+\|\nabla\eta\|_\infty^2\Big)
\|\1_{\eta w\not= 0}\|_{r;s}^2.  
\end{align}
\el
\begin{proof}
Since $\alpha\in[0,1]$, by relation \eqref{NH19}, one sees that
$$
1-\frac{1}{p}=\frac{\alpha}{d}+\alpha\left(\frac{r+2}{2r}-\frac{1}{d}\right)+\frac{2(1-\alpha)}{r}.
$$
Thus by  H\"older's inequality and Gagliardo-Nirenberge's inequality \eqref{EQ1}, we have
\begin{align}\label{AS1}
\begin{split}
\<g, h\>&=\left|\int_{\mR^d}(\mI-\Delta)^{-\alpha/2} g\cdot (\mI-\Delta)^{\alpha/2}  h\right|\\
&\leq\|g\|_{-\alpha,p}\|h\|_{\alpha,p/(p-1)}\lesssim \|g\|_{-\alpha,p}\|h\|^{1-\alpha}_{r/2}\|\nabla h\|^{\alpha}_{2r/(r+2)}.
\end{split}
\end{align}
\\
(i) Since $\chi_2|_{Q_2}\equiv 1$ and $\eta|_{Q^c_2}=0$, by 
\eqref{AS1} with $g=c\chi_2$ and $h=\eta^2 w^2$, we have
\begin{align*}
\<c,\eta^2w^2\>=\<c\chi_2,\eta^2w^2\>
&\lesssim  \|c\chi_2\|_{-\alpha,p}\|\eta^2 w^2\|^{1-\alpha}_{r/2}\|\nabla(\eta^2 w^2)\|^{\alpha}_{2r/(r+2)}\\
&\lesssim  \|c\chi_2\|_{-\alpha,p}\|\eta w\|^{2(1-\alpha)}_{r}\|(\eta w)\nabla(\eta w)\|^{\alpha}_{2r/(r+2)}\\
&\lesssim  \|c\chi_2\|_{-\alpha,p}\|\eta w\|^{2-\alpha}_{r}\|\nabla(\eta w)\|^{\alpha}_2,
\end{align*}
where we drop the time variable $t$ and  the last step is due to H\"older's inequality.
Integrating both sides in the time variable, and 
due to  $\frac{2-\a}{s}+\frac{\a}{2}+\frac{1}{q}=1$, by H\"older's inequality again we get 
$$
|\<\!\<c, \eta^2 w^2\>\!\>|
\lesssim \|c\chi_2\|_{-\alpha,p;q}\|\eta w\|^{2-\alpha}_{r;s}\|\nabla(\eta w)\|^{\alpha}_{2;2},
$$
which gives \eqref{NJ52} by Young's inequality.
\ \medskip\\
(ii) By \eqref{AS1} with $g=b\chi_2$ and $h=\nabla\eta^2 w^2$, we have
$$
\<b, \nabla\eta^2 w^2\>=\<b\chi_2, \nabla\eta^2 w^2\>
\lesssim  \|b\chi_2\|_{-\alpha,p}\|\nabla\eta^2 w^2\|^{1-\alpha}_{r/2}\|\nabla(\nabla\eta^2 w^2)\|^{\alpha}_{2r/(r+2)}.
$$
Notice that 
$$
\|\nabla\eta^2 w^2\|_{r/2}\leq 2\|\nabla\eta\|_\infty\|w\1_{\eta\not=0}\|_r^2,
$$
and  by H\"older's inequality,
\begin{align*}
\|\nabla(\nabla\eta^2 w^2)\|_{2r/(r+2)}
&\leq\|\nabla^2\eta^2 w^2\|_{2r/(r+2)}+\|\nabla\eta^2 \nabla w^2\|_{2r/(r+2)}\\
&\leq\|\nabla^2\eta^2\|_{2r/(r-2)}\|w\1_{\eta\not=0}\|^2_r+4\|w\nabla\eta\|_r\cdot\|\eta\nabla w\|_2.
\end{align*}
Hence,
\begin{align*}
\<b, \nabla\eta^2 w^2\>&\lesssim\|b\chi_2\|_{-\alpha,p}
\|\nabla\eta\|^{1-\alpha}_\infty\|w\1_{\eta\not=0}\|^{2(1-\alpha)}_r
\Big(\|\nabla^2\eta^2\|_{2r/(r-2)}^\alpha\| w\1_{\eta\not=0}\|^{2\alpha}_r+\|w\nabla\eta\|^{\alpha}_r\cdot\|\eta\nabla w\|^{\alpha}_2\Big)\\
&\lesssim \|b\chi_2\|_{-\alpha,p}\Big(\big(\|\nabla\eta\|^{1-\alpha}_\infty\|\nabla^2\eta\|^{\alpha}_\infty+\|\nabla\eta\|^{1+\alpha}_\infty\big)\|w\1_{\eta\not=0}\|^{2}_r
+\|\nabla\eta\|_\infty\|w\1_{\eta\not=0}\|^{2-\alpha}_r\|\eta\nabla w\|^{\alpha}_2\Big)\\
&\lesssim \|b\chi_2\|_{-\alpha,p}\Big(\big(1+\|\nabla^2\eta\|_\infty+\|\nabla\eta\|^2_\infty\big)\|w\1_{\eta\not=0}\|^{2}_r
+\|\nabla\eta\|_\infty\|w\1_{\eta\not=0}\|^{2-\alpha}_r\|\nabla(\eta w)\|^{\alpha}_2\Big),
\end{align*}
and by H\"older's inequality and due to $\frac{1}{q}+\frac{2-\a}{s}+\frac{\a}{2}=1$, $s\geq \frac{2q}{q-1}$, 
\begin{align*}
|\<\!\<b, \nabla\eta^2 w^2\>\!\>|
&\lesssim \|b\chi_2\|_{-\alpha,p;q}(1+\|\nabla^2\eta\|_\infty+\|\nabla\eta\|^2_\infty)\|w\1_{\eta\not=0}\|^{2}_{r;2q/(q-1)}\\
&\quad+\|b\chi_2\|_{-\alpha,p;q}\|\nabla\eta\|_\infty\|w\1_{\eta\not=0}\|^{2-\alpha}_{r;s}\|\nabla(\eta w)\|^{\alpha}_{2;2},
\end{align*}
The desired estimate \eqref{NJ51} follows by Young's inequality and $2q/(q-1)\leq s$.
\ \medskip\\
(iii)
By \eqref{AS1} with $g=f\eta$ and $h=\eta w$, we have
\begin{align*}
\<f\eta,\eta w\>\lesssim \|f\eta\|_{-\alpha,p}\|\eta w\|^{1-\alpha}_{r/2}\|\nabla(\eta w)\|^\alpha_{2r/(r+2)}.
\end{align*}
Since $\nabla(\eta w)=\nabla(\eta w)^+-\nabla (\eta w)^-=0$ on $\{\eta w=0\}$ (cf. \cite[Lemma 7.6]{Gi-Tr}), we have 
\begin{align}\label{LK}
\nabla(\eta w)=\nabla(\eta w) \1_{\eta w\not=0},\ \ a.s.
\end{align}
Thus, by H\"older's inequality, we further have
\begin{align*}
|\<f\eta,\eta w\>|&\lesssim\|f\eta\|_{-\alpha,p}\|\eta w\|^{1-\alpha}_{r/2}\|\nabla(\eta w)\|^{\alpha}_2\|\1_{\eta w\not=0}\|_r^\a\\
&\lesssim\|f\eta\|_{-\alpha,p}\|\eta w\|^{1-\alpha}_{r}\|\nabla(\eta w)\|^{\alpha}_2\|\1_{\eta w\not=0}\|_r,
\end{align*}
and  due to $\frac{1}{q}+\frac{1-\alpha}{s}+\frac{\alpha}{2}+\frac{1}{s}=1$, $\frac{d}{r}+\frac{2}{s}>\frac{d}{2}$,
\begin{align}\label{EB1}
|\<\!\<f, \eta^2 w\>\!\>|&\lesssim\|f\eta\|_{-\alpha,p;q}\|\eta w\|_{r;s}^{1-\alpha}\|\nabla(\eta w)\|^{\alpha}_{2;2} \|\1_{\eta w\not=0}\|_{r;s}
\stackrel{\eqref{EB}}{\lesssim}\|f\eta\|_{-\alpha,p;q}\|\eta w\|_{\sV}\|\1_{\eta w\not=0}\|_{r;s}.
\end{align}
Notice that by $\eta=\chi_2\eta$ and \eqref{EW2} below,
$$
\|f\eta\|_{-\alpha,p;q}\lesssim\|f\chi_2\|_{-\alpha,p;q}\|\eta\|_{1,\infty}\lesssim \|f\chi_2\|_{-\alpha,p;q}(1+\|\nabla\eta\|_{\infty}).
$$
Substituting this into \eqref{EB1} and by Young's inequality, we obtain \eqref{NJ55}.
\end{proof}

\subsection{Local energy estimate}
Throughout this paper we shall always assume 
$$
b\in\mL^2_{2, loc},\ \ f\in\sD',
$$
and consider the following PDE in $\mR^{d+1}$:
\begin{align}\label{PDE0}
\p_t u=\Delta u+b\cdot \nabla u+f.
\end{align}
\bd
A function $u\in\sV_{loc}\cap\mL^\infty_{loc}$ is called a weak solution
(subsolution or supersolution)  
of PDE \eqref{PDE0} with coefficients $(b,f)$
if for any nonnegative smooth function $\varphi\in C^\infty_c(\mR^{d+1})$ and almost all $t\in \mR$, 
\begin{align}\label{Def0}
\begin{aligned}
\<\p_tu,\varphi\>=(\leq \mbox{or }\geq)-\<\nabla u,\nabla\varphi\>+\<b\cdot\nabla u,\varphi\>+\<f,\varphi\>.
\end{aligned}
\end{align}

\ed

Now we prove the following local energy estimate.
\bl[Energy estimate]
Suppose that for some $(\alpha_i,p_i,q_i)\in\sI_d$,  $i=1,2,3$,
$$
b\in \mH^{-\alpha_1,p_1}_{q_1,loc},\ \ -\div b\leq \varTheta_b \in \mH^{-\alpha_2,p_2}_{q_2,loc},\ \ f\in \mH^{-\alpha_3,p_3}_{q_3,loc}.
$$
Let $(r_i,s_i)\in[2,\infty]$ be defined by \eqref{NH19} and $\kappa\geq 0$.
For any weak subsolution $u\in\sV_{loc}\cap\mL^\infty_{loc}$ 
of PDE \eqref{PDE0}, there is a constant $C>0$ depending only on $d,\alpha_i,p_i,q_i, i=1,2,3$ and
$$
\|b\chi_2\|_{-\alpha_1,p_1;q_1},\ \ \| \varTheta_b\chi_2\|_{-\alpha_2,p_2; q_2},
$$ 
where $\chi_2$ is defined by \eqref{CH},
such that for $w:=(u-\kappa)^+$ and any $\eta\in C^\infty_c(Q_2)$ and $t\geq 0$,
\begin{align}\label{Es1}
\|\eta w\cI_t\|_{\sV}\leq C 
\Xi_\eta^{1/2}\left(\|w\1_{\eta\not=0}\cI_t\|_{r_1;s_2}+\|w\eta\cI_t\|_{r_2;s_2}
+\|f\chi_2\cI_t\|_{-\alpha_3,p_3;q_3}\|\1_{w\eta\not=0}\|_{r_3;s_3}\right),
\end{align}
where $\cI_t(\cdot):=\1_{(-\infty,t]}(\cdot)$, and
\begin{align}\label{Xi}
\Xi_\eta:=1+\|\p_t\eta\|_\infty+
\|\nabla\eta\|_\infty^2+\|\nabla^2\eta\|_\infty.
\end{align}
\el
\begin{proof}
By taking the Steklov mean of $u$, without loss of generality we may assume $\p_t u \in \mL^2_{2,loc}$.
By \eqref{Def0} and smoothing approximation, for any nonnegative $\varphi\in\sV_{loc}\cap\mL^\infty_{loc}$ with compact support in $Q_2$, we have
$$
\<\p_tu,\varphi\>\leq-\<\nabla u,\nabla\varphi\>+\<b\cdot\nabla u,\varphi\>+\<f,\varphi\>.
$$
Let $\eta\in C^\infty_c(Q_2)$ and $w:=(u-\kappa)^+$ for some $\kappa\geq 0$. 
Taking the test function $\varphi=\eta^2 w\in\sV_{loc}\cap\mL^\infty_{loc}$ and integrating in time variable from $-\infty$ to $t$, 
by the integration by parts formula, we have
\begin{align}\label{HA1}
\begin{split}
\int_{\mR^d}(\eta w)^2(t)&\leq 2\int_{\Gamma_t}\!\p_t\eta^2 w^2- 2\int^t_{-\infty}\!\!\<\nabla u,\nabla(\eta^2w)\>
+ 2\int^t_{-\infty}\!\!\<b\cdot \nabla u, \eta^2w\>+ 2\int^t_{-\infty}\!\!\<f,\eta^2w\>,
\end{split}
\end{align}
where $\Gamma_t:=(-\infty,t]\times\mR^d$ and we have used that 
$$
\<\p_t u,\eta^2 w\>=\<\p_t w,\eta^2 w\>=\frac{1}{2}\left[\p_t\int_{\mR^d}(\eta w)^2-\int_{\mR^d}\p_t\eta^2w^2\right].
$$
Noticing that
\begin{align}\label{HA2}
\nabla u\cdot\nabla w=|\nabla w|^2,\ (\nabla u) w=(\nabla w)w=\nabla w^2/2,
\end{align}
by  the integration by parts formula again, we have
\begin{align*}
2\<\nabla u,\nabla(\eta^2w)\>
&=2 \int_{\mR^d}\eta^2|\nabla w|^2- \int_{\mR^d} w^2\cdot\Delta\eta^2,
\end{align*}
and
\begin{align*}
2\<b\cdot \nabla u, \eta^2w\>=\<b, \nabla w^2\eta^2\>=- \<\div b, \eta^2w^2\>- \<b,\nabla \eta^2w^2\>.
\end{align*}
Therefore, by $-\div b\leq\Theta_b$ and smoothing approximation for $w$, we further have
\begin{align*}
&\int_{\mR^d}(\eta w)^2(t)+2 \int_{\Gamma_t}\eta^2|\nabla w|^2\leq  \int_{\Gamma_t}(|2\p_t\eta^2+\Delta\eta^2|)w^2\\
&\qquad+\left|\int_{-\infty}^t \<\varTheta_b, \eta^2w^2\>\right| +\left|\int_{-\infty}^t \<b, \nabla \eta^2w^2\>\right|+2\left|\int_{-\infty}^t  \<f, \eta^2w\>\right|, 
\end{align*}
which yields by definition \eqref{NK1} that
\begin{align*}
\|\eta w\cI_t\|^2_{\sV}&\leq \int_{\Gamma_t}(|2\p_t\eta^2+\Delta\eta^2|+4|\nabla\eta^2|)w^2
+\sup_{s\leq t}\left|\int_{-\infty}^s \<\varTheta_b, \eta^2w^2\>\right| \\
&\quad+\sup_{s\leq t}\left|\int_{-\infty}^s \<b, \nabla \eta^2w^2\>\right| +2\sup_{s\leq t}\left|\int_{-\infty}^s  \<f, \eta^2w\>\right|=:\sum_{i=1}^4 I_i.
\end{align*}
For $I_1$,  notice that $r_1, s_1\geq 2$, we have 
$$
I_1\leq\Big(4\|\p_t\eta\|_\infty+2\|\nabla\eta\|^2_\infty+\|\Delta\eta\|_\infty+4\|\nabla\eta\|_\infty\Big) \|w\1_{\eta\not=0}\cI_t\|_{2;2}^2 
\leq C \Xi_\eta\|w\1_{\eta\not=0}\cI_t\|_{r_1;s_1}^2. 
$$
For $I_2, I_3$ and $I_4$, by \eqref{NJ52}, \eqref{NJ51} and \eqref{NJ55}, we have
\begin{align*}
I_2\leq \eps\|\eta w\cI_t\|^2_{\sV}+C_\eps\| \varTheta_b\chi_2\|^2_{-\alpha_2,p_2;q_2} \|w\eta\cI_t\|_{r_2;s_2}^2,
\end{align*}
and
\begin{align*}
I_3\leq \eps\|\eta w\cI_t\|^2_{\sV}+C_\eps\Big(1+\|b\chi_2\|^2_{-\alpha_1,p_1;q_1}\Big)\Big(1+\|\nabla^2\eta\|_\infty+\|\nabla\eta\|_\infty^2\Big)
\|w\1_{\eta\not=0}\cI_t\|_{r_1;s_1}^2,
\end{align*}
\begin{align*}
I_4\leq \eps\|\eta w\cI_t\|^2_{\sV}+C_\eps\|f\chi_2\cI_t\|^2_{-\alpha_3,p_3;q_3}(1+\|\nabla\eta\|_\infty^2) \|\1_{\eta w \not= 0}\|_{r_3;s_3}^2.  
\end{align*}
Combining the above calculations and letting $\eps$ be small enough, we obtain
\begin{align*}
\|\eta w\cI_t\|^2_{\sV}\lesssim  
\Xi_\eta\|w\1_{\eta\not=0}\cI_t\|_{r_1;s_1}^2+\|w\eta\cI_t\|_{r_2;s_2}^2
+(1+\|\nabla\eta\|_\infty^2)\|f\chi_2\cI_t\|^2_{-\alpha_3,p_3;q_3}  
\ \|\1_{\eta w \not= 0}\|_{r_3;s_3}^2,   
\end{align*}
where  $\Xi_\eta$ is define by \eqref{Xi}.
From this, we derive \eqref{Es1}. 
\end{proof}

\br\label{Re23}
If $\alpha_1=0$ and $\frac{d}{p_1}+\frac{2}{q_1}=1$ or $b(t,x)=b(x)\in L^d_{loc}(\mR^d)$, 
then we can remove the assumption on the divergence of $b$. In fact, in this case,
we can give a direct treatment for the term $b$ in \eqref{HA1} as follows: For any $\eps>0$, 
let 
$$
b_\eps(t,x):=b(t,\cdot)*\rho_\eps(x),\ \ \widetilde b_\eps(t,x):=b(t,x)-b_\eps(t,x).
$$ 
Let $\frac{1}{r_1}+\frac{1}{p_1}=\frac{1}{s_1}+\frac{1}{q_1}=\frac{1}{2}$, which satisfy $\frac{d}{r_1}+\frac{2}{s_1}=\frac{d}{2}$ due to $\frac{d}{p_1}+\frac{2}{q_1}=1$. 
Since $\chi_2\eta=\eta$, by \eqref{HA2}, H\"older's inequality and Lemma \ref{Le11}, we have
\begin{align*}
\int_{\Gamma_t}\left|(b\cdot \nabla u) \eta^2w\right|&\leq\int_{\Gamma_t}\left|(\widetilde b_\eps\cdot \nabla u) \eta^2w\right|+
\int_{\Gamma_t}\left|(b_\eps\cdot \nabla u) \eta^2w\right|\\
&\leq \|\widetilde b_\eps\chi_2\|_{p_1;q_1}\|\eta\nabla w\cI_t\|_{2;2}\|\eta w\cI_t\|_{r_1;s_1}\\
&\quad+\|b_\eps\chi_2\|_{\infty;\infty}\|\eta\nabla w\cI_t\|_{2;2}\|\eta w\cI_t\|_{2;2}\\
&\leq c_\eps\|\eta w\cI_t\|^2_\sV+C_\eps\|\1_{\eta\not=0} w\cI_t\|^2_{2;2},
\end{align*}
where $\lim_{\eps\to 0}c_\eps=0$ and $\lim_{\eps\to 0}C_\eps=\infty$.
Using this estimate to replace the corresponding estimate about $b$ and taking $\eps$ small enough, we still have
\eqref{Es1}.
Here the reason that for $p=d$ we assume $b$ being time-independent is that in general
$$
\lim_{\eps\to 0}\|\widetilde b_\eps\chi_2\|_{d;\infty}\not=0\ \mbox{for $b\in L^\infty_{loc}(L^d_{loc})$, but }\ 
\lim_{\eps\to 0}\|\widetilde b_\eps\chi_2\|_{d}=0\ \mbox{ for $b\in L^d_{loc}$}.
$$
\er
\subsection{Maximum principle}
The following De Giorgi's iteration lemma is well known \cite{Ha-Li}.
\bl\label{Le15}
Let $(a_n)_{n\in\mN}$ be a sequence of nonnegative numbers. Suppose that for some $C_0,\lambda>1$ and $\eps>0$,
$$
a_{n+1}\leq C_0\lambda^n a_n^{1+\eps},\ \ n=1,2,\cdots.
$$
If $a_1\leq C_0^{-1/\eps}\lambda^{-1/\eps^2}$, then
$$
\lim_{n\to\infty} a_n=0.
$$
\el

Now we can show the following local maximum principle for PDE \eqref{PDE0}.
\bt[Local maximum estimate]
Suppose that for some $(\alpha_i,p_i,q_i)\in\sI_d$,  $i=1,2,3$,
$$
b\in\mH^{-\alpha_1,p_1}_{q_1;loc},\ \ -\div b\leq\varTheta_b\in\mH^{-\alpha_2,p_2}_{q_2;loc},\ \ f\in\mH^{-\alpha_3,p_3}_{q_3;loc}.
$$
For any weak subsolution $u\in\sV_{loc}\cap\mL^\infty_{loc}$ of PDE \eqref{PDE0}, 
there is a constant $C>0$ depending only on $d,\alpha_i,p_i,q_i, i=1,2,3$ and the quantities
$$
\|b\chi_2\|_{-\alpha_1,p_1;q_1},\ \ \|\varTheta_b\chi_2\|_{-\alpha_2,p_2; q_2},
$$ 
where $\chi_2$ is defined by \eqref{CH},  such that
\begin{align}\label{JK2}
\|u^+\1_{Q_1}\|_{\infty}\leq C\left(\|u^+\chi_2\|_{\sV}+\|f\chi_2\|_{-\alpha_3,p_3;q_3}\right).
\end{align}
\et
\begin{proof}
Let $\kappa>0$, which will be determined below. For $n\in\mN$, define
$$
t_n:=4\cdot(4^{-1}+3\cdot 4^{-n}),\ \ \lambda_n:=1+2^{1-n},\ \ \kappa_n:=\kappa\left(1-2^{1-n}\right)
$$
and
$$
\Gamma_n:=(-t_n,t_n)\times B_{\lambda_n}\downarrow [-1,1]\times \bar B_1=\bar Q_1.
$$
Let $\zeta^{\rm t}_n\in C^\infty_c(-2,2)$ be a time-cutoff function so that for some $C>0$ and any $n\in\mN$,
$$
\zeta^{\rm t}_n|_{(-t_{n+1}, t_{n+1})}=1,\ \ \zeta^{\rm t}_n|_{(-t_{n}, t_{n})^c}=0,\ \  |\p_t\zeta^{\rm t}_n|\leq C4^n.
$$
Let $\zeta^{\rm x}_n\in C^\infty_c(B_2)$ be a spatial-cutoff function so that for some $C>0$ and any $n\in\mN$,
$$
\zeta^{\rm x}_n|_{B_{\lambda_{n+1}}}=1,\ \ \zeta^{\rm x}_n|_{B_{\lambda_n}^c}=0,\ \  |\nabla^j\zeta^{\rm x}_n|\leq C2^{jn},\ \ j=1,2.
$$
 Now let us define
$$
\eta_n(t,x):=\zeta^{\rm t}_n(t)\cdot\zeta^{\rm x}_n(x).
$$
Let $\Xi_{\eta_n}$ be defined by \eqref{Xi}. It is easy to see that for some $C>0$ and all $n\in\mN$,
$$
\eta_n|_{\Gamma_{n+1}}=1,\ \ \eta_n|_{\Gamma_n^c}=0,\ \
\Xi_{\eta_n}\leq C 4^n.
$$
Let
$$
w_n:=(u-\kappa_n)^+.
$$
Notice that
$$
w_n|_{w_{n+1}\not=0}=(u-\kappa_{n+1}+\kappa_{n+1}-\kappa_n)^+|_{w_{n+1}\not=0}\geq \kappa_{n+1}-\kappa_{n}=\kappa 2^{-n}.
$$
For $i=1,2,3$, due to $\eta_n|_{\Gamma^c_n}=0$, we have
$$
\ell^{(i)}_n:=\|w_n\1_{\Gamma_n}\|_{r_i;s_i}
\geq\|w_n\1_{\eta_nw_{n+1}\not=0}\|_{r_i;s_i}\geq \kappa 2^{-n}\|\1_{\eta_{n}w_{n+1}\not=0}\|_{r_i;s_i},
$$
which means that
\begin{align}\label{LK1}
\|\1_{\eta_{n}w_{n+1}\not=0}\|_{r_i;s_i}\leq 2^{n}\ell^{(i)}_n /\kappa.
\end{align}
Since $\frac{2}{r_i}+\frac{d}{s_i}>\frac{d}{2}$, we can choose $\gamma_{i},\beta_{i}>r, \theta_{i},\tau_{i}>s$ so that
$$
\frac{1}{\gamma_i}+\frac{1}{\beta_i}=\frac{1}{r_i},\quad \frac{1}{\theta_i}+\frac{1}{\tau_i}=\frac{1}{s_i},\quad \frac{d}{\gamma_i}+\frac{2}{\theta_i}\geq \frac{d}{2}.
$$
Thus, by $\eta_n|_{\Gamma_{n+1}}=1$, H\"older's inequality, Lemmas \ref{Le11} and \ref{Le12}, we have
\begin{align}\label{LK3}
\begin{split}
\ell^{(i)}_{n+1}&=\|w_{n+1}\1_{\Gamma_{n+1}}\|_{r_i;s_i}\leq \|\eta_nw_{n+1}\|_{r_i;s_i}\\
&\leq\|\eta_nw_{n+1}\|_{\gamma_i;\theta_i}\|\1_{\eta_nw_{n+1}\not=0}\|_{\beta_i;\tau_i}\\
&\leq\|\eta_nw_{n+1}\|_{\gamma_i;\theta_i}\|\1_{\eta_nw_{n+1}\not=0}\|^{(s_i/\tau_i)\wedge(r_i/\beta_i)}_{r_{i};s_{i}}\\
&\leq C\|\eta_nw_{n+1}\|_{\sV}\cdot(2^{n}\ell^{(i)}_n /\kappa)^{(s_i/\tau_i)\wedge(r_i/\beta_i)}.
\end{split}
\end{align}
Notice $\Gamma_1=Q_2$. By \eqref{Es1} with $\eta=\eta_n$ and $w=w_{n+1}$, for  $\kappa\geq\|f\chi_2\|_{-\alpha_3,p_3;q_3}$, we obtain
\begin{align}\label{LK33}
\begin{split}
\|\eta_nw_{n+1}\|_{\sV}&\lesssim 2^n\left(\|w_{n+1}\1_{\Gamma_n}\|_{r_1;s_1}+\|w_{n+1}\eta_n\|_{r_2;s_2}+\|f\chi_2\|_{-\a_3,p_3;q_3}\|\1_{\eta_nw_{n+1}\not=0}\|_{r_3;s_3}\right)\\
&\stackrel{\eqref{LK1}}{\lesssim} C 2^n(\ell^{(1)}_n+\ell^{(2)}_n+2^n\ell^{(3)}_n)\lesssim  4^n(\ell^{(1)}_n+\ell^{(2)}_n+\ell^{(3)}_n).
\end{split}
\end{align}
Now we put
$$
a_n:=\Big(\ell^{(1)}_n+\ell^{(2)}_n+\ell^{(3)}_n\Big)/\kappa.
$$
By \eqref{LK3} and \eqref{LK33}, we obtain that for some $C_0,\eps>0$ and $\lambda>1$,
\begin{align*}
a_{n+1}\leq C4^n a_n\sum_{i=1}^3(2^na_n)^{(s_i/\tau_i)\wedge(r_i/\beta_i)}\leq C_0\lambda^n a_n^{1+\eps},
\end{align*}
provided  $\kappa\geq\|f\chi_2\|_{-\alpha_3,p_3;q_3}$.
Notice that by $\chi_2|_{\Gamma_1}=1$ and Lemma \ref{Le11},
\begin{align*}
a_1&\leq \frac{1}{\kappa}\sum_{i=1}^3\|u^+\1_{\Gamma_1}\|_{r_i;s_i}\leq \frac{1}{\kappa}\sum_{i=1}^3\|u^+\chi_2\|_{r_i;s_i}\leq \frac{C_1}{\kappa}\|u^+\chi_2\|_{\sV}.
\end{align*}
If $\kappa\geq (C_1C_0^{1/\eps}\lambda^{1/\eps^2}\|u^+\chi_2\|_{\sV})\vee\|f\chi_2\|_{-\alpha_3,p_3;q_3}$ 
so that $a_1\leq C_0^{-1/\eps}\lambda^{-1/\eps^2}$, then by Fatou's lemma and Lemma \ref{Le15},
\begin{align*}
\|(u-\kappa)^+\1_{Q_1}\|_{r_1;s_1}\leq \liminf_{n\to\infty}\|w_n\1_{\Gamma_n}\|_{r_1;s_1}
=\liminf_{n\to\infty}\ell^{(1)}_n\leq\kappa\cdot\limsup_{n\to\infty}a_n=0,
\end{align*}
which implies that for Lebesgue almost all $(t,x)\in\mR\times\mR^d$,
$$
(u^+\1_{Q_1})(t,x)\leq C_1C_0^{1/\eps}\lambda^{1/\eps^2}\|u^+\chi_2\|_{\sV}\vee\|f\chi_2\|_{-\alpha_3,p_3;q_3}.
$$
The proof is complete.
\end{proof}
\br\label{Re24}
If $\alpha_1=0$ and $\frac{d}{p_1}+\frac{2}{q_1}=1$ or $b(t,x)=b(x)\in L^d_{loc}(\mR^d)$, 
then by Remark  \ref{Re23}, we can drop the condition on the divergence of $b$.
\er
Now we aim to prove the following crucial result.
\bt\label{TH22}
(Global maximum estimate) Suppose that for some $(\alpha_i,p_i,q_i)\in\sI_d$,  $i=1,2,3$,
\begin{align}\label{Con3}
b\in\widetilde\mH^{-\alpha_1,p_1}_{q_1},\ \ -\div b\leq\varTheta_b\in\widetilde\mH^{-\alpha_2,p_2}_{q_2},\ \ f\in\widetilde\mH^{-\alpha_3,p_3}_{q_3}.
\end{align}
Let  $u\in\sV_{loc}\cap\mL^\infty_{loc}$ be a weak solution of PDE \eqref{PDE0}  with initial value $u(0)=0$.
For any $T>0$, there exists a constant $C>0$ depending only on $T, d,\alpha_i,p_i,q_i$ and the quantity
$$
\kappa:=\nor b\nor_{-\alpha_1,p_1;q_1}+\nor \varTheta_b\nor_{-\alpha_1,p_1;q_1}
$$ 
such that 
\begin{align}\label{Es161}
\|u\|_{L^\infty([0,T]\times\mR^d)}+\nor u \1_{[0,T]}\nor_{\sV}\leq C\nor f\1_{[0,T]}\nor_{-\alpha_3,p_3;q_3}.
\end{align}
\et
\begin{proof}
Without loss of generality, we assume $T=1$ and 
$$
u(t,x)=f(t,x)\equiv 0, \ \  \forall t\leq 0. 
$$
Let $\chi_1$ be as in \eqref{CH} and define for $z\in\mR^d$,
$$
\eta_z(t,x):=\chi_1(t,x-z).
$$
By translation and  
\eqref{Es1} with $\eta=\eta_z$ and $w=u^+, u^-$, there is a constant 
$C>0$ depending only on $T, d$, $\alpha_i,p_i,q_i$, $\nor b\nor_{-\alpha_1,p_1;q_1}$, $\nor (\div b)^-\nor_{-\alpha_2,p_2;q_2}$ such that for all $t\in[0,1]$,
$$
\|\eta_z u \cI_t\|_{\sV}\leq C\Big(\|u\1_{\eta_z\not=0}\cI_t\|_{r_1;s_1}+\|\eta_z u\cI_t\|_{r_2;s_2}
+\|f\chi^{0,z}_2\cI_t\|_{-\alpha_3,p_3;q_3}\Big),
$$
where $\chi^{0,z}_2$ is the same as in \eqref{CH}.
Taking supremum in $z\in\mR^d$ for both sides, we obtain
\begin{align}\label{Es131}
\sup_z\|\eta_z u \cI_t\|_{\sV}\leq C\Big(\sup_z\|u\1_{\eta_z\not=0}\cI_t\|_{r_1;s_1}+\sup_z\|\eta_z u\cI_t\|_{r_2;s_2}
+\sup_z\| f\chi^{0,z}_2\cI_t\|_{-\alpha_3,p_3;q_3}\Big).
\end{align}
Since for each $z\in\mR^d$, there are at most $N$-points $z_1,\cdots, z_N\in\mR^d$ such that
$$
B_2(z)\subset\cup_{j=1}^NB_1(z_j),
$$
where $N=N(d)$, we have for $t\in[0,1]$,
\begin{align}\label{Es141}
\|u\1_{\eta_z\not=0}\cI_t\|_{r_1;s_1}\leq \|u\1_{B_2(z)}\cI_t\|_{r_1;s_1}
\leq\sum_{j=1}^N\|u\1_{B_1(z_j)}\cI_t\|_{r_1;s_1}\leq N\sup_z\|\eta_z u\cI_t\|_{r_1;s_1},
\end{align}
where the last step is due to $\eta_{z_j}|_{[0,t]\times B_1(z_j)}=1$. Hence, by \eqref{Es131}, \eqref{Es141} and \eqref{GT1} in appendix,
\begin{align}\label{Es101}
\sup_z\|\eta_z u \cI_t\|_{\sV}\leq C\Big(\sup_z\|\eta_z u\cI_t\|_{r_1;s_1}+\sup_z\|\eta_z u\cI_t\|_{r_2;s_2}
+\nor f\cI_t\nor_{-\alpha_3,p_3;q_3}\Big).
\end{align}
Let 
$$
\theta_i:=\frac{d}{2}-\frac{d}{r_i},\ \ s_i':=\frac{2(1-\theta_i)s_i}{2-s_i\theta_i}.
$$ 
Since $\frac{d}{r_i}+\frac{2}{s_i}>\frac{d}{2}$, by Lemma \ref{Le11}, we have for any $\eps\in(0,1)$,
\begin{align}\label{Es121}
\|\eta_z u\cI_t\|_{r_i;s_i}\leq \eps\|\nabla (\eta_z u)\cI_t\|_{2;2}+C_\eps\|\eta_z u\cI_t\|_{2;s'_i}.
\end{align}
Combining \eqref{Es101} and \eqref{Es121}, we arrive at
\begin{align*}
&\sup_{z}\|\eta_z u \cI_t\|_{2;\infty}+\sup_{z}\|\nabla(\eta_z u) \cI_t\|_{2;2}\leq 2\sup_{z}\|\eta_z u \cI_t\|_{\sV}\\
&\quad\leq \eps\sup_z\|\nabla (\eta_z u)\cI_t\|_{2;2}+C_\eps\sup_z\|\eta_z u\cI_t\|_{2;(s'_1\vee s'_2)}
+C\nor f\cI_t\nor_{-\alpha_3,p_3;q_3}.
\end{align*}
By choosing $\eps$ small enough, we obtain
\begin{align}\label{Es151}
\sup_{z}\|\eta_z u \cI_t\|_{2;\infty}+\sup_{z}\|\nabla(\eta_z u) \cI_t\|_{2;2}\leq
C\sup_z\|\eta_z u\cI_t\|_{2;(s'_1\vee s'_2)}+C\nor f\cI_t\nor_{-\alpha_3,p_3;q_3}.
\end{align}
Since $s_1', s_2'<\infty$ and $u(t)\equiv 0$ for $t\leq 0$,  the above inequality implies that for any $t\in[0,1]$,
$$
\sup_{z}\|(\eta_z u)(t)\|^{s'_1\vee s'_2}_{2}\leq C\sup_z\int^t_0\|(\eta_z u)(s)\|^{s'_1\vee s'_2}_{2}\dif s+C\nor f\cI_t\nor^{s'_1\vee s'_2}_{-\alpha_3,p_3;q_3}.
$$
By Gronwall's inequality we obtain 
$$
\sup_{z}\sup_{t\in[0,1]}\|(\eta_z u)(t)\|_{2}\leq C\nor f\1_{[0,1]}\nor_{-\alpha_3,p_3;q_3},
$$
which together with \eqref{Es151} yields 
$$
\nor u\1_{[0,1]}\nor_{\sV}\lesssim
\sup_{z}\|\eta_z u\1_{[0,1]}\|_{2;\infty}+\sup_{z}\|\nabla(\eta_z u)\1_{[0,1]}\|_{2;2}\lesssim\nor f\1_{[0,1]}\nor_{-\alpha_3,p_3;q_3}.
$$
Finally,  by \eqref{JK2} and \eqref{Es161}, we also have
\begin{align*}
\|u\|_{L^\infty([0,1]\times\mR^d)}&\leq \sup_{z}\|(u^++u^{-})\1_{[0,1]\times B^z_1}\|_{\infty}
\lesssim \nor u\1_{[0,1]}\nor_{\sV}+\nor f\1_{[0,1]}\nor_{-\alpha_3,p_3;q_3} \lesssim \nor f\1_{[0,1]}\nor_{-\alpha_3,p_3;q_3},
\end{align*}
where $B^z_1:=\{x: |x-z|<1\}$.
The proof is complete.
\end{proof}
\subsection{Existence-uniqueness and stability}
In this subsection we prove the existence-uniqueness and stability of weak solutions for PDE \eqref{PDE0} by using the apriori estimate \eqref{Es161}.
For $T>0$ and a function $f$ in $\mR^{d+1}$, we denote
$$
f^T:=f\1_{[0,T]}, \ \widetilde\sV_T:=\{f: \nor f^T\nor_{\sV}<\infty\},\ \mL^\infty_T:=\{f: \|f^T\|_\infty<\infty\}.
$$ 

\bt\label{TH23}
(Existence-uniqueness) Under \eqref{Con3}, for any $T>0$,
there exists a unique weak solution $u\in\widetilde\sV_T\cap\mL^\infty_T$ to PDE \eqref{PDE0} with initial value $u(0)=0$.
\et
\begin{proof}
First of all, the uniqueness is a direct consequence of \eqref{Es161}. We prove the existence by weak convergence method.
Let $b_n(t,x):=b(t,\cdot)*\rho_n(x)$ and $f_n(t,x):=f(t,\cdot)*\rho_n(x)$. By (ii) of Proposition \ref{Pr41} in Appendix, we have
$$
b_n\in L^{q_1}_{loc}(\mR_+; C^\infty_b(\mR^d)), \ \ f_n\in L^{q_3}_{loc}(\mR_+; C^\infty_b(\mR^d)). 
$$
and
\begin{align}\label{GK1}
\sup_n\left(\nor b_n\nor_{-\alpha_1, p_1;q_1}+\nor (\div b_n)^-\nor_{-\alpha_2,p_2;q_2}+\nor f_n\nor_{-\alpha_3,p_3;q_3}\right)<\infty.
\end{align}
It is well known that the following PDE has a unique smooth solution $u_n\in C(\mR_+; C^\infty_b(\mR^d))$:
\begin{align}\label{FQ1}
\p_t u_n=\Delta u_n+b_n\cdot\nabla u_n+f_n=0.
\end{align}
By \eqref{GK1} and Theorem \ref{TH22}, we have
\begin{align}\label{EM1}
\sup_n\left(\|u^T_n\|_{\infty}+\nor u^T_n\nor_{\sV}\right)<\infty,\ \ \forall T>0.
\end{align}
Hence, by the fact that every bounded subset of 
$\widetilde \mH^{1,2}_2$ 
is relatively compact, there is a subsequence $n_k$ and $\bar u\in\cap_{T>0}(\widetilde\sV_T\cap\mL^\infty_T)$ 
such that for any $\varphi\in C^\infty_c(\mR^{d+1})$ and $g\in \mH^{-1,2}_{2;loc}$, 
\begin{align}\label{EM2}
\lim_{k\to\infty}\<\!\< u_{n_k},  g \varphi\>\!\>=\<\!\<\bar u, g \varphi\>\!\>.
\end{align}
by taking weak limits for equation \eqref{FQ1}, one finds that $\bar u$ is a weak solution of PDE \eqref{PDE0}.
Indeed, it suffices to prove that for any $\varphi\in C_c^\infty(\mR^{d+1})$, 
\begin{align}\label{EM3}
&\lim_{k\to\infty}\<\!\< b_{n_k} \cdot\nabla u_{n_k},\varphi\>\!\>=\<\!\< b \cdot\nabla \bar u,\varphi\>\!\>,\quad
\lim_{k\to\infty}\<\!\< f_{n_k},\varphi\>\!\>=\<\!\< f,\varphi\>\!\>.
\end{align}
Let the support of $\varphi$ be contained in $Q_R$ for some $R>0$. 
Since $b\in\mL^2_{2;loc}$, by \eqref{EM1} and H\"older's inequality, we have
for some $C>0$ independent of $k$,
\begin{align*}
\<\!\<(b_{n_k}-b) \cdot\nabla u_{n_k},\varphi\>\!\>&=\<\!\<\chi_R(b_{n_k}-b)\cdot\nabla u_{n_k},\varphi\>\!\>
\leq\|\nabla \varphi\|_\infty\|(b_{n_k}-b)\chi_R\|_{2;2}\|\chi_R\nabla u_{n_k}\|_{2;2}\\
&\leq C \|(b_{n_k}-b)\chi_R\|_{2;2}\to 0\ \mbox{ as\ \  $k\to\infty$.}
\end{align*}
On the other hand, since $\div(b\varphi)\in\mH^{-1,2}_2$ has compact support, by \eqref{EM2} we also have
$$
\lim_{k\to\infty}\<\!\<b\cdot\nabla(u_{n_k}-\bar u),\varphi\>\!\>=
\lim_{k\to\infty}\<\!\<u_{n_k}-\bar u,\div(b\varphi)\>\!\>=0.
$$
Thus we obtain the first limit in \eqref{EM3}. The second limit in \eqref{EM3} is direct.
\end{proof}

\bt(Stability) \label{TH24}
Let $(p_i,q_i)\in[2,\infty)$ with $\frac{d}{p_i}+\frac{2}{q_i}<2$, where $i=1,2,3$. For any $n\in\mN\cup\{\infty\}=:\mN_\infty$, 
let $b_n, f_n\in\sD'$ satisfy
\begin{align}\label{GK11}
\sup_{n\in\mN_\infty}\left(\nor b_n\nor_{p_1;q_1}+\nor (\div b_n)^-\nor_{p_2;q_2}+\nor f_n\nor_{p_3;q_3}\right)<\infty.
\end{align}
For $n\in\mN_\infty$, let $u_n\in\widetilde\sV\cap\mL^\infty$ be the unique weak solutions of PDE \eqref{PDE0} associated with coefficients $(b_n,f_n)$ with initial value $u(0)=0$.
Assume that for any $\varphi\in C_c(\mR^{d+1})$,
\begin{align}\label{Lim0}
\lim_{n\to\infty}\left(\|(b_n-b)\varphi\|_{p_1;q_1}+\|\div (b_n-b)\varphi\|_{p_2;q_2}+\|(f_n-f)\varphi\|_{p_3;q_3}\right)=0.
\end{align}
Then it holds that for Lebesgue almost all $(t,x)\in\mR_+\times\mR^{d}$,
\begin{align}\label{Lim}
\lim_{n\to\infty}u_n(t,x)= u_\infty(t,x).
\end{align}
\et
\begin{proof}
Notice that equation
\begin{align*}
\p_t u_n=\Delta u_n+b_n\cdot\nabla u_n+f_n=\Delta u_n+\div (b_nu_n)-(\div b_n) u_n+f_n
\end{align*}
holds in the distributional sense (see \eqref{Def0}).
Letting $r:=\frac{2p_1}{p_1+2}$ and $s:=\frac{2q_1}{q_1+2}$, by Proposition \ref{Pr41} in Appendix, we have
\begin{align*}
\nor(\p_t u_n)\1_{[0,T]}\nor_{-1,r;s}&\leq\nor\Delta u^T_n+\div (b^T_nu^T_n)-(\div b^T_n) u^T_n+f^T_n\nor_{-1,r;s}\\
&\lesssim\nor u^T_n\nor_{1,r;s}+\nor b^T_nu^T_n\nor_{r;s}+\nor(\div b^T_n) u^T_n\nor_{-1,r;s}+\nor f^T_n\nor_{-1, r;s}\\
&\lesssim\nor u^T_n\nor_{1,2;2}+\nor b^T_n\nor_{r;s}\|u^T_n\|_\infty+\nor \div b^T_n\nor_{-1,p_1;q_1}\nor u^T_n\nor_{1,2;2}+\nor f^T_n\nor_{r;s}\\
&\lesssim\nor u^T_n\nor_{1,2;2}+\nor b^T_n\nor_{p_1;q_1}\left(\|u^T_n\|_\infty+\nor u^T_n\nor_{1,2;2}\right)+\nor f^T_n\nor_{p_3;q_3}.
\end{align*}
By \eqref{GK11} and Theorem \ref{TH22}, we get
$$
\sup_n\Big(\|u_n\|_{\mL^\infty_T}+\nor u_n\nor_{\sV_T}+\nor(\p_t u_n)\1_{[0,T]}\nor_{-1,r;s}\Big)<\infty.
$$
Thus by Aubin-Lions' lemma (cf. \cite{Si}), there is a subsequence $n_k$ and $\bar u\in \cap_{T>0}(\widetilde\sV_T\cap\mL^\infty_T)$ such that
$$
\lim_{k\to\infty}\|u_{n_k}-\bar u\|_{L^2([0,T]\times B_m)}=0,\ \ \forall T>0, m\in\mN.
$$
By selecting a subsubsequence $n'_k$, it holds that for Lebesgue almost all $(t,x)\in\mR_+\times\mR^d$,
\begin{align}\label{Lim1}
u_{n'_k}(t,x)\to \bar u(t,x),\ \ k\to\infty.
\end{align}
As in showing \eqref{EM3}, one can show that $\bar u$ is a weak solution of PDE \eqref{PDE0}.
By the uniqueness, $\bar u=u$, and by a contradiction method, the whole sequence converges almost everywhere. 
\end{proof}

\section{Proof of Theorem \ref{TH11}}

Below we always assume that for some $p_i,q_i\in[2,\infty)$ with $\frac{d}{p_i}+\frac{2}{q_i}<2$, $i=1,2$,
$$
\kappa:=\nor b\nor_{p_1;q_1}+\nor (\div b)^-\nor_{p_2;q_2}<\infty.
$$
Let $b_n(t,x)=b(t,\cdot)*\rho_n(x)$ be the mollifying approximation of $b(t,\cdot)$. 
By (ii) of Proposition \ref{Pr41} in Appendix, we have
\begin{align}\label{HF1}
\sup_n\left(\nor b_n\nor_{p_1;q_1}+\nor (\div b_n)^-\nor_{p_2;q_2}\right)\leq C\kappa,
\end{align}
and
$$
b_n\in L^{q_1}_{loc}(\mR_+;C^\infty_b(\mR^d)).
$$
For $(s,x)\in\mR_+\times\mR^d$, consider the following SDE:
\begin{align}\label{SDE9}
\dif X^n_{s,t}=b_n(t,X^n_{s,t})\dif t+\sqrt{2}\dif W_t,\ \ X^n_{s,s}=x,\ \ t\geq s,
\end{align}
where $W$ is a $d$-dimensional standard Brownian motion on some complete filtered probability space
$(\Omega,\sF,(\sF_t)_{t\geq 0},\bP)$.
It is well known that there is a unique strong solution $X^n_{s,t}(x)$ to the above SDE (cf. \cite{Ik-Wa}).

\medskip
Now we are in the position to prove our main result and we divide our proof into three parts. 

\subsection{Existence of martingale solutions}
First of all, we prove the following crucial  estimate of Krylov's type.
\bl
For any $(\alpha,p,q)\in\sI_d$, there are constants $\theta=\theta(\alpha,p,q)>0$ 
and $C>0$ depending on $\kappa,d$, $\alpha,p,q$, $p_i,q_i$ such that for any $f\in C^\infty_c(\mR^{d+1})$
and $0\leq s\leq  t_0<t_1<\infty$ with $t_0-t_1\leq 1$,
\begin{align}\label{Kry}
\sup_n\sup_{x\in\mR^d}\bE\left(\int^{t_1}_{t_0} f(t,X_{s,t}^n(x))\dif t\Bigg|\sF_{t_0}\right)\leq C(t_1-t_0)^\theta\nor f\nor_{-\alpha,p;q}.
\end{align}
In particular, we have the following Khasminskii's estimate: for any $\lambda\in\mR$,
\begin{align}\label{Kas}
\sup_n\sup_{x\in\mR^d}\bE\exp\left\{\lambda\int^{s+1}_s |f(t,X_{s,t}^n(x))|\dif t\right\}\leq C=C(\lambda,\kappa,\nor f\nor_{-\alpha,p;q}).
\end{align}
\el
\begin{proof}
Fix $0\leq s\leq t_0<t_1<\infty$ with $t_0-t_1\leq 1$ and $f\in C^\infty_c(\mR^{d+1})$. Let $u_n$ be the smooth solution of the following backward PDE:
\begin{align}\label{EG1}
\p_t u_n+\Delta u_n+b_n\cdot\nabla u_n+f=0,\ u_n(t_1,\cdot)=0.
\end{align}
By It\^o's formula we have
$$
u_n(t_1,X^n_{s,t_1})=u_n(t_0,X^n_{s,t_0})+\int^{t_1}_{t_0}
(\p_t u_n+\Delta u_n+b_n\cdot\nabla u_n)(t,X^n_{s,t})\dif t+\sqrt{2}\int^{t_1}_{t_0}\nabla u_n(t,X^n_{s,t})\dif W_t.
$$
By \eqref{EG1} and taking conditional expectation with respect to $\sF_{t_0}$, we obtain
\begin{align}\label{EM8}
\bE\left(\int^{t_1}_{t_0}f(t,X^n_{s,t})\dif t\Bigg|\sF_{t_0}\right)=\bE \Big(u_n(t_0,X^n_{s,t_0})|\sF_{t_0}\Big)\leq\|u_n(t_0)\|_\infty.
\end{align}
Since $\frac{d}{p}+\frac{2}{q}<2-\alpha$, we can choose $q'<q$ so that $\frac{d}{p}+\frac{2}{q'}<2-\alpha$.
Thus by Theorem \ref{TH22} and H\"older's inequality, there is constant $C=C(\kappa,d,\alpha,p,q,p_i,q_i)>0$ such that
$$
\bE\left(\int^{t_1}_{t_0}f(t,X^n_{s,t})\dif t\Bigg|\sF_{t_0}\right)
\leq C\nor f\1_{[t_0,t_1]}\nor_{-\alpha,p;q'}\leq C(t_1-t_0)^{1-\frac{q'}{q}}\nor f\nor_{-\alpha,p;q}.
$$
Thus we obtain \eqref{Kry}. As for \eqref{Kas}, it is a direct consequence of \eqref{Kry} and \cite[Lemma 1.1]{Po} 
(or see \cite{Zh0}).
\end{proof}

\bl\label{Le32}
For any $T>0$, there is a constant $C>0$ such that for any $f\in L^1(\mR^d)$ and $n\in\mN$,
$$
\|\cT^n_{s,t}f\|_1\leq C\|f\|_1,\ \ \forall 0\leq s<t\leq s+T,
$$
where $\cT^n_{s,t}f(x):=\bE f(X^n_{s,t}(x))$.
Moreover, if $(\div b)^{-}\equiv 0$, then the above $C$ can be $1$.
\el
\begin{proof}
Let $Y^n_{s,t}:=Y^n_{s,t}(x)$
be the inverse flow of $x\mapsto X^n_{s,t}(x)$. 
Notice that $s\mapsto Y^n_{s,t}$ solves the following backward SDE:
$$
Y^n_{s,t}=x-\int_s^t b_n(r,Y^n_{r,t})\dif r+\sqrt{2} (W_s-W_t),\ 0\leq s\leq t.
$$
Letting $J^n_{s,t}:=J^n_{s,t}(x):=\nabla Y^n_{s,t}(x)$ be the Jacobian matrix, we have
$$
\p_s J^n_{s,t}=\nabla b_n(s,Y^n_{s,t})J^n_{s,t}\Rightarrow \p_s \det(J^n_{s,t})=\div b_n(s,Y^n_{s,t})\det(J^n_{s,t}).
$$
Hence,
$$
\det (J^n_{s,t})=\exp\left\{-\int^t_s\div b_n(r,Y^n_{r,t})\dif r\right\}\leq\exp\left\{\int^t_s(\div b_n)^{-}(r,Y^n_{r,t})\dif r\right\}.
$$
Thus, by Khasminskii's estimate \eqref{Kas} and \eqref{HF1}, we have
$$
\sup_n\sup_{x\in\mR^d}\bE \det(J^n_{s,t}(x))<\infty.
$$
Now by the change of variables, for any nonnagative $f\in L^1(\mR^d)$, we have
$$
\|\cT^n_{s,t}f\|_1=\bE\left( \int_{\mR^d}f(X^n_{s,t}(x))\dif x\right)=\bE\left( \int_{\mR^d}f(x)\det(J^n_{s,t}(x))\dif x\right)\leq C\|f\|_1.
$$
Moreover, if $(\div b)^{-}\equiv 0$, then $\det (J^n_{s,t})\leq 1$ and the above $C\equiv1$.
\end{proof}

\bl\label{Le33}
For each $(s,x)\in\mR_+\times\mR^d$, let $\mP^n_{s,x}$ be the law of $X^n_{s,\cdot}(x)$ in $\mC$. Then
$(\mP^n_{s,x})_{n\in\mN}$ is tight.
\el
\begin{proof}
Fix $(s,x)\in\mR_+\times\mR^d$ and $T>s$. Let $\tau\geq s$ be any stopping time less than $T$. Notice that
$$
X^n_{s,\tau+\delta}(x)-X^n_{s,\tau}(x)=\int^{\tau+\delta}_\tau b_n(t,X^n_{s,t}(x))\dif t+\sqrt{2}(W_{\tau+\delta}-W_\tau),\ \ \delta>0.
$$
By the strong Markov property and \eqref{Kry} with $\alpha=0$, we have
\begin{align*}
\bE|X^n_{s,\tau+\delta}(x)-X^n_{s,\tau}(x)|&
\leq\sup_{t\in[s,T], y\in\mR^d}\bE\left(\int^{t+\delta}_t |b_n(r,X^n_{t,r}(y))|\dif r\right)+\sqrt{2}\delta^{1/2}\\
&\leq C\delta^\theta\nor b_n\nor_{p;q}+\sqrt{2}\delta^{1/2}\leq C\delta^\theta\nor b\nor_{p;q}+\sqrt{2}\delta^{1/2},
\end{align*}
where $C$ is independent of $n$.
Thus by \cite[Lemma 2.7]{Zh-Zh1}, we obtain
$$
\sup_n\sup_{(s,x)\in[0,T]\times\mR^d}\bE\left(\sup_{t\in[s,T]}|X^n_{s,t+\delta}(x)-X^n_{s,t}(x)|^{1/2}\right)\leq C\left(\delta^{\theta/2}\nor b\nor_{p;q}^{1/2}+\delta^{1/4}\right).
$$
From this, by Chebyshev's inequality, we derive that for any $\eps>0$,
$$
\lim_{\delta\to 0}\sup_n\sup_{(s,x)\in[0,T]\times\mR^d}\bP\left(\sup_{t\in[s,T]}|X^n_{s,t+\delta}(x)-X^n_{s,t}(x)|>\eps\right)=0.
$$
Hence, by \cite[Theorem 1.3.2]{St-Va}, the law of $X^n_{s,\cdot}(x)$ is tight in $\mC$.
\end{proof}
Now we can show the existence of martingale solutions.
\bl\label{Le34}
Any accumulation point $\mP_{s,x}$ of $(\mP^n_{s,x})_{n\in\mN}$ belongs to $\sM_{s,x}^b$. Moreover, for any $(\alpha,p,q)\in\sI_d$,
there are $\theta=\theta(\alpha,p,q)>0$ and constant $C>0$ such that for any $f\in C^\infty_c(\mR^{d+1})$ and $0\leq s\leq  t_0<t_1<\infty$ with $t_1-t_0\leq 1$,
\begin{align}\label{KL}
\sup_{x\in\mR^d}\mE^{\mP_{s,x}}\left(\int^{t_1}_{t_0} f(t,\omega_t)\dif t\Big|\cB_{t_0}\right)\leq C(t_1-t_0)^\theta\nor f\nor_{-\alpha,p;q}.
\end{align}
\el
\begin{proof}
Let $(\alpha,p,q)\in\sI_d$. By \eqref{Kry}, 
there are $\theta=\theta(\alpha,p,q)>0$ and constant $C>0$ such that for each $n\in\mN$, $x\in\mR^d$, and any $f\in C^\infty_c(\mR^{d+1})$, 
$0\leq s\leq  t_0<t_1<\infty$ with $t_1-t_0\leq 1$, and
$G\in C_b(\mC)$ being $\cB_{t_0}$-measurable, 
\begin{align}\label{HF2}
\mE^{\mP^n_{s,x}}\left(\int^{t_1}_{t_0} f(t,\omega_t)\dif t\cdot G_{t_0}\right)\leq C(t_1-t_0)^\theta\nor f\nor_{-\alpha,p;q}\mE (G_{t_0}).
\end{align}
Let $\mP_{s,x}$ be any accumulation point of $(\mP^n_{s,x})_{n\in\mN}$, that is, for some subsequence $n_k$,
$$
\mP^{n_k}_{s,x}\mbox{ weakly converges to $\mP_{s,x}$ as $k\to\infty$.}
$$
By taking weak limits for \eqref{HF2} and a standard monotone class method, we obtain \eqref{KL}.
In order to prove $\mP_{s,x}\in\sM_{s,x}^b$, it suffices to prove that for any $t_1>t_0\geq s$ and $f\in C^2_c(\mR^d)$,
$$
\mE^{\mP_{s,x}}(M^f_{t_1}|\cB_{t_0})=M^f_{t_0},\ \ \ \mP_{s,x}-a.s.,
$$
where 
$$
M^f_t:=f(\omega_t)-f(x)-\int^t_s(\Delta+b\cdot\nabla) f(r,\omega_r)\dif r.
$$
By the standard monotone class method, it is enough to show that for any $G\in C_b(\mC)$ being $\cB_{t_0}$-measurable, 
$$
\mE^{\mP_{s,x}}\Big(M^f_{t_1}\cdot G_{t_0}\Big)=\mE^{\mP_{s,x}}\Big(M^f_{t_0}\cdot G_{t_0}\Big).
$$
Note that for each $n\in\mN$,
$$
\mP^n_{s,x}\in\sM^{b_n}_{s,x}\Rightarrow\mE^{\mP^n_{s,x}}\Big(M^f_{t_1}\cdot G_{t_0}\Big)
=\mE^{\mP^n_{s,x}}\Big(M^f_{t_0}\cdot G_{t_0}\Big).
$$
We want to take weak limits, where the key point is to show
\begin{align}\label{FU2}
\lim_{k\to\infty}\mE^{\mP^{n_k}_{s,x}}
\left(\int^{t_1}_s(b\cdot\nabla f)(r,\omega_r)\dif r\cdot G_{t_0}(\omega)\right)=\mE^{\mP_{s,x}}
\left(\int^{t_1}_s(b\cdot\nabla f)(r,\omega_r)\dif r\cdot G_{t_0}(\omega)\right).
\end{align}
Assume that supp$(f)\subset Q_R$. 
By \eqref{Kry} with $\alpha=0$ and \eqref{GT3} in Appendix, we have
\begin{align}\label{FU1}
\begin{split}
&\sup_n\mE^{\mP^n_{s,x}}\left(\int^{t_1}_s((b_m-b)\cdot\nabla f)(r,\omega_r)\dif r\cdot G_{t_0}(\omega)\right)\\
&\quad\leq\|G_{t_0}\|_\infty\|\nabla f\|_\infty\sup_n\mE^{\mP^n_{s,x}}
\left(\int^{t_1}_s|(b_m-b)\chi_R|(r,\omega_r)\dif r\right)\\
&\quad\lesssim\|G_{t_0}\|_\infty\|\nabla f\|_\infty\nor(b_m-b)\chi_R\nor_{p_1;q_1}\to 0,\ m\to\infty,
\end{split}
\end{align}
where $\chi_R|_{Q_R}\equiv1 $ has compact support.
On the other hand, for fixed $m\in\mN$, since
$$
\omega\mapsto \int^{t_1}_s(b_m\cdot\nabla f)(r,\omega_r)\dif r\cdot G_{t_0}(\omega)\in C_b(\mC),
$$
we also have
$$
\lim_{k\to\infty}\mE^{\mP^{n_k}_{s,x}}\left(\int^{t_1}_s(b_m\cdot\nabla f)(r,\omega_r)\dif r\cdot G_{t_0}(\omega)\right)
=\mE^{\mP_{s,x}}\left(\int^{t_1}_s(b_m\cdot\nabla f)(r,\omega_r)\dif r\cdot G_{t_0}(\omega)\right),
$$
which together with \eqref{FU1} yields \eqref{FU2}.
The proof is complete.
\end{proof}

\subsection{Weak convergence of $\mP^n_{s,x}$}
In this subsection we show that for Lebesgue almost all $(s,x)$,
the accumulation point of $(\mP^n_{s,x})_{n\in\mN}$ is unique, which in turn implies that
$$
\mP^n_{s,x}\mbox{ weakly  converges to }\mP_{s,x}\in\sM^b_{s,x}\mbox{ as $n\to\infty$}.
$$
For fixed $T>0$ and $f\in \mL^\infty_T$, by Theorem \ref{TH23}, 
there is a unique weak solution $u=u_{T,f}\in\widetilde\sV_T\cap \mL^\infty_T$ 
to the following backward PDE:
\begin{align}\label{PDE1}
\p_t u+\Delta u+b\cdot\nabla u+f=0,\ u(t,\cdot)|_{t\geq T}=0.
\end{align}
Let $\mQ\subset\mR$ be the set of all rational numbers and $\sG_0$ be a countable dense subset of $C^\infty_c(\mR^{d})$.
For $m\in\mN$, we recursively define a countable set $\sG_m$ as follows:
$$
\sG_m:=\Big\{g=f u_{T,h}\in \mL^\infty_T: T\in\mQ,\ f\in\sG_0,\ h\in\sG_{m-1}\Big\}.
$$
Clearly,
$$
\sA:=\cup_{m=0}^\infty\sG_m\subset\mL^\infty\mbox{ is a countable set.}
$$

\bl\label{Le35}
For $T>0$, $f\in\mL^\infty_T$ and $n\in\mN$, if we define
\begin{align}\label{UU}
u^n_{T,f}(s,x)=\mE^{\mP^n_{s,x}}\left(\int^T_sf(t,\omega_t)\dif t\right),
\end{align}
then $u^n_{T,f}\in\widetilde\sV_T\cap\mL^\infty_T$ uniquely solves PDE \eqref{PDE1} with $b=b_n$. Moreover,
there is a Lebesgue null set $\cN\subset \mR_+\times\mR^d$ such that  for  all $(s,x)\in \cN^c$, $f\in\sA$ and $s\leq T\in\mQ$,
\begin{align}\label{LIM}
\lim_{n\to\infty}u^{n}_{T,f}(s,x)=u_{T,f}(s,x).
\end{align}
\el
\begin{proof}
For $m\in\mN$, let $f_m(t,x):=f(t,\cdot)*\rho_m(x)$ and
$$
u^{n,m}_{T,f}(s,x)=\mE^{\mP^n_{s,x}}\left(\int^T_sf_m(t,\omega_t)\dif t\right),\ s\in[0,T], \ x\in\mR^d.
$$
It is well known that $u^{n,m}_{T,f}$ solves PDE \eqref{PDE1} with $b=b_n$ and $f=f_m$.
By Theorem \ref{TH24}, for Lebesgue almost all $(s,x)$, we have
\begin{align}\label{UU9}
u^{n,m}_{T,f}(s,x)\to u^n_{T,f}(s,x),
\end{align}
where $u^n_{T,f}\in\widetilde\sV_T\cap\mL^\infty_T$ is the unique weak solution of of PDE \eqref{PDE1} with $b=b_n$. 
On the other hand, by Krylov's estimate \eqref{Kry}, for each $s\leq T$ and $x\in\mR^d$, we have
$$
\lim_{m\to\infty}\mE^{\mP^n_{s,x}}\left(\int^T_sf_m(t,\omega_t)\dif t\right)=\mE^{\mP^n_{s,x}}\left(\int^T_sf(t,\omega_t)\dif t\right),
$$
which together with \eqref{UU9} gives \eqref{UU}. 
For fixed $T\in\mQ$ and $f\in\sG_0$, by Theorem \ref{TH24} again, there is a Lebesgue null set $\cN_{T,f}\subset \mR_+\times\mR^d$ such that 
\eqref{LIM} holds for all $(s,x)\in \cN^c_{T,f}$. Finally, we just need to take
$$
\cN:=\cup_{T\in\mQ}\cup_{f\in\sG_0}\cN_{T,f}.
$$
The proof is complete.
\end{proof}

\bl\label{Le36}
Let $\cN$ be as in Lemma \ref{Le35}.
For fixed $(s,x)\in\cN^c$ and any two accumulation points $\mP^{(1)}_{s,x}$ and $\mP^{(2)}_{s,x}$ of $(\mP^n_{s,x})_{n\in\mN}$,
\begin{align}\label{GH2}
\mP^{(1)}_{s,x}=\mP^{(2)}_{s,x}.
\end{align}
\el
\begin{proof}
Fix $(s,x)\in\cN^c$. For $s\leq T\in\mQ$ and $f\in \sG_0$, 
by \eqref{LIM} and taking weak limits for \eqref{UU} along different subsequences for $\mP^{(i)}_{s,x}$, $i=1,2$, one finds that 
$$
u_{T,f}(s,x)=\mE^{\mP^{(i)}_{s,x}}\left(\int^T_sf(\omega_t)\dif t\right),\ \ i=1,2,
$$
which implies that for all $s\leq T\in\mQ$ and $f\in \sG_0$,
$$
\mE^{\mP^{(1)}_{s,x}}\left(\int^T_sf(\omega_t)\dif t\right)=\mE^{\mP^{(2)}_{s,x}}\left(\int^T_sf(\omega_t)\dif t\right).
$$
In particular, for all $T\geq s$ and $f\in\sG_0$,
$$
\int^T_s\mE^{\mP^{(1)}_{s,x}}f(\omega_t)\dif t=\int^T_s\mE^{\mP^{(2)}_{s,x}}f(\omega_t)\dif t
\Rightarrow \mE^{\mP^{(1)}_{s,x}} f(\omega_T)=\mE^{\mP^{(2)}_{s,x}} f(\omega_T).
$$
{\it Claim:} Let $g_m(t,x)\in\mL^\infty$ be uniformly bounded and converge to $g(t,x)$ for Lebesgue almost all $(t,x)$. 
For any $(s,x)\in\cN^c$ and $T>s$, it holds that for each $x\in\mR^d$,
\begin{align}\label{KP2}
\lim_{m\to\infty}\sup_n\mE^{\mP^n_{s,x}}\left(\int^T_s|g_m-g|(t,\omega_t)\dif t\right)=0.
\end{align}
{\it Proof of Claim:} For $R>0$, define
$$
\tau_R:=\inf\{t\geq s: |\omega_t|\geq R\}.
$$
By \eqref{Kry} with $(\alpha,p,q)=(0,d,4)$ and the dominated convergence theorem, we have
\begin{align}\label{LL9}
\begin{split}
\lim_{m\to\infty}\sup_n\mE^{\mP^n_{s,x}}\left(\int^{T\wedge\tau_R}_s|g_m-g|(t,\omega_t)\dif t\right)
\leq C\lim_{m\to\infty}\big\|(g_m-g)\1_{[s,T]\times B_R}\big\|_{d;4}=0.
\end{split}
\end{align}
On the other hand, by SDE \eqref{SDE9} and \eqref{Kry} again, we also have
$$
\bE\left(\sup_{t\in[s,T]}|X^n_{s,t}|\right)\lesssim|x|+1+\bE\left(\int^{T}_0|b_n(t,X_{s,t}^n)| \dif t\right)\leq C,
$$
where $C$ is independent of $n$. Hence,
$$
\lim_{R\to\infty}\sup_n\mP^n_{s,x}(\tau_R<T)
=\lim_{R\to\infty}\sup_n\bP\left(\sup_{t\in[s,T]}|X^n_{s,t}(x)|>R\right)\leq\lim_{R\to\infty}\sup_n\bE\left(\sup_{t\in[s,T]}|X^n_{s,t}|\right)/R=0,
$$
which together with \eqref{LL9} yields the claim.

\medskip

Next let $s\leq T_1<T_2$ be two rational numbers and $f_1,f_2\in\sG_0$. 
Let $(\mP^{n_k}_{s,x})_{k\in\mN}$ be a subsequence so that $(\mP^{n_k}_{s,x})_{k\in\mN}$ weakly converges to $\mP^{(1)}_{s,x}$. 
By the Markov property, we have for $f_1,f_2\in\sG_0$
\begin{align*}
&\mE^{\mP^{(1)}_{s,x}}\left(\int^{T_1}_sf_1(\omega_{t_1})\left(\int^{T_2}_{t_1}f_2(\omega_{t_2})\dif t_2\right)\dif t_1\right)\\
&=\lim_{k\to\infty}\mE^{\mP^{n_k}_{s,x}}\left(\int^{T_1}_sf_1(\omega_{t_1})\left(\int^{T_2}_{t_1}f_2(\omega_{t_2})\dif t_2\right)\dif t_1\right)\\
&=\lim_{k\to\infty}\mE^{\mP^{n_k}_{s,x}}\left(\int^{T_1}_sf_1(\omega_{t_1})\mE^{\mP^{n_k}_{t_1,\omega_{t_1}}}\left(\int^{T_2}_{t_1}f_2(\omega_s)\dif s\right)\dif t_1\right)\\
&\stackrel{\eqref{UU}}{=}\lim_{k\to\infty}\mE^{\mP^{n_k}_{s,x}}\left(\int^{T_1}_sf_1(\omega_{t_1}) u^{n_k}_{T_2,f_2}(t_1,\omega_{t_1})\dif t_1\right)\\
&=\lim_{k\to\infty}\mE^{\mP^{n_k}_{s,x}}\left(\int^{T_1}_sf_1(\omega_{t_1}) u_{T_2,f_2}(t_1,\omega_{t_1})\dif t_1\right),
\end{align*}
where the last step is due to \eqref{LIM} and the above Claim.
Notice that
$$
g(s,x):=f_1(x)u_{T_2,f_2}(s,x)\in\sA.
$$
Hence,
$$
\lim_{k\to\infty}\mE^{\mP^{n_k}_{s,x}}\left(\int^{T_1}_sf_1(\omega_{t_1}) u_{T_2,f_2}(t_1,\omega_{t_1})\dif t_1\right)
\stackrel{\eqref{UU}}{=}\lim_{k\to\infty}u^{n_k}_{T_1,g}(s,x)\stackrel{\eqref{LIM}}{=}u_{T_1,g}(s,x).
$$
Since the right hand side does not depend on the choice of the subsequence $n_k$, we finally obtain
that for any rational numbers $s\leq T_1<T_2$  and $f_1,f_2\in\sG_0$,
$$
\mE^{\mP^{(1)}_{s,x}}\left(\int^{T_1}_sf_1(\omega_{t_1})\left(\int^{T_2}_{t_1}f_2(\omega_{t_2})\dif t_2\right)\dif t_1\right)
=\mE^{\mP^{(2)}_{s,x}}\left(\int^{T_1}_sf_1(\omega_{t_1})\left(\int^{T_2}_{t_1}f_2(\omega_{t_2})\dif t_2\right)\dif t_1\right).
$$
From this, as above we derive that for all $f_1,f_2\in\sG_0$ and $T_2>T_1\geq s$,
$$
\mE^{\mP^{(1)}_{s,x}}\left(f_1(\omega_{T_1})f_2(\omega_{T_2})\right)
=\mE^{\mP^{(2)}_{s,x}}\left(f_1(\omega_{T_1})f_2(\omega_{T_2})\right).
$$
Similarly, we can prove that for any $T_m>\cdots >T_1\geq s$ and $f_1,\cdots, f_m\in\sG_0$,
$$
\mE^{\mP^{(1)}_{s,x}}\left(f_1(\omega_{T_1})\cdots f_m(\omega_{T_m})\right)
=\mE^{\mP^{(2)}_{s,x}}\left(f_1(\omega_{T_1})\cdots f_m(\omega_{T_m})\right).
$$
Thus we obtain \eqref{GH2}.
\end{proof}
\subsection{Almost surely Markov property}
Let $\cN$ be as in Lemma \ref{Le35}.
We fix $(s,x)\in\cN^c$ so that 
\begin{align}\label{JG1}
\mP^n_{s,x}\mbox{ weakly converges to $\mP_{s,x}$ as $n\to\infty$}.
\end{align}
Recalling that $\sG_0$ is a countable dense subset of $C^\infty_c(\mR^{d})$, to show \eqref{JH2}, it suffices to prove the following claim: 
\\
\\
{\it Claim 1:} 
For fixed $t_1\in(s,\infty)\cap\mQ$ and $f\in\sG_0$, there is a Lebesgue-null set $I^{t_1,f}_{s,x}\subset(s,t_1)$
so that for all $t_0\in(s,t_1)\setminus I^{t_1,f}_{s,x}$,
\begin{align}\label{CON}
\mE^{\mP_{s,x}}\left( f(\omega_{t_1})|\cB_{t_0}\right)=\mE^{\mP_{t_0,\omega_{t_0}}}\left( f(\omega_{t_1})\right),\ \ \mP_{s,x}-a.s.
\end{align}
Indeed, if this is proven, then we can take
$$
I_{s,x}:=\cup_{s<t_1\in \mQ}\cup_{f\in\sG_0}I^{t_1,f}_{s,x}.
$$
Thus for any $t_0\in(s,\infty)\setminus I_{s,x}$, and all $t_0<t_1\in\mQ$ and $f\in\sG_0$,
$$
\mE^{\mP_{s,x}}\left( f(\omega_{t_1})|\cB_{t_0}\right)=\mE^{\mP_{t_0,\omega_{t_0}}}\left( f(\omega_{t_1})\right),\ \ \mP_{s,x}-a.s.
$$
By a standard approximation argument, the above equality also holds for all $t_1>t_0$ and $f\in C_c(\mR^d)$.
\\
\\
Furthermore, to prove Claim 1, it suffices to prove the following claim:
\\
\\
{\it Claim 2:} Let $t_1\in(s,\infty)\cap\mQ$ and $f\in\sG_0$. For fixed $m\in\mN$, $s_1,\cdots,s_m\in(s,t_1)\cap\mQ$ and $g_1,\dots, g_m\in \sG_0$,
there exists a null set $I:=I^{s_1,\cdots,s_m}_{g_1,\cdots,g_m}\subset[s_m,t_1]$ so that for all $t_0\in[s_m,t_1]\setminus I$,
\begin{align}\label{KP1}
\mE^{\mP_{s,x}}\left(g_1(\omega_{s_1})\cdots g_m(\omega_{s_m})f(\omega_{t_1})\right)
=\mE^{\mP_{s,x}}\left(g_1(\omega_{s_1})\cdots g_m(\omega_{s_m})\mE^{\mP_{t_0,\omega_{t_0}}}(f(\omega_{t_1}))\right).
\end{align}
Indeed, if this is proven, then we can define
$$
I^{t_1,f}_{s,x}:=\cup_{m\in\mN}\cup_{s_1,\cdots,s_m\in (s,t_1)\cap\mQ}\cup_{g_1,\cdots,g_m\in\sG_0}I^{s_1,\cdots,s_m}_{g_1,\cdots,g_m}\subset(s,t_1).
$$
Thus for any $t_0\in(s,t_1)\setminus I^{t_1,f}_{s,x}$, \eqref{KP1} holds for all $m$ and 
$s_1,\cdots, s_m\in(s,t_0]\cap\mQ$, $g_1,\cdots,g_m\in\sG_0$. By a standard monotone class argument, we obtain \eqref{CON}
for $t_0\in(s,t_1)\setminus I^{t_1,f}_{s,x}$ from \eqref{KP1}.
\\
\\
{\it Proof of Claim 2:} For simplicity of notations, we shall write
$$
G_{s_m}(\omega):=g_1(\omega_{s_1})\cdots g_m(\omega_{s_m}).
$$
By the Lebesgue differential theorem, we only need to prove that for any $t_0\in[s_m,t_1]$,
\begin{align}\label{JG7}
\mE^{\mP_{s,x}}\left(G_{s_m}(\omega)f(\omega_{t_1})\right)
=\frac{1}{t_1-t_0}\int^{t_1}_{t_0}\mE^{\mP_{s,x}}\left(G_{s_m}(\omega)\mE^{\mP_{r,\omega_{r}}}(f(\omega_{t_1}))\right)\dif r.
\end{align}
Clearly, by the Markov property of $(\mP^n_{s,x})_{(s,x)\in\mR_+\times\mR^d}$, we have
\begin{align}\label{JG4}
\mE^{\mP^n_{s,x}}\left(G_{s_m}(\omega)f(\omega_{t_1})\right)
=\frac{1}{t_1-t_0}\int^{t_1}_{t_0}\mE^{\mP^n_{s,x}}\left(G_{s_m}(\omega)\mE^{\mP^n_{r,\omega_{r}}}(f(\omega_{t_1}))\right)\dif r.
\end{align}
By \eqref{JG1} we have
\begin{align}\label{JG5}
\lim_{n\to\infty}\mE^{\mP^n_{s,x}}\left(G_{s_m}(\omega)f(\omega_{t_1})\right)=\mE^{\mP_{s,x}}\left(G_{s_m}(\omega)f(\omega_{t_1})\right).
\end{align}
Define
$$
H_n(r,y):=\mE^{\mP^n_{r,y}} f(\omega_{t_1})=\bE f(X^n_{r,t_1}(y)).
$$
Since by \eqref{JG1}, $H_n(r,y)\to H(r,y)$ for Lebesgue almost all $r,y$, by \eqref{KP2}, we have
$$
\lim_{n\to\infty}\sup_k\mE^{\mP^k_{s,x}}\left(\int^{t_1}_{t_0} |H_n(r,\omega_r)-H(r,\omega_r)|\dif r\right)=0.
$$
On the other hand, for fixed $n$, since $y\mapsto H_n(r,y)$ is continuous, we also have
$$
\lim_{k\to\infty}\int^{t_1}_{t_0}\mE^{\mP^k_{s,x}}\left(G_{s_m}(\omega)H_n(r,\omega_{r})\right)\dif r
=\int^{t_1}_{t_0}\mE^{\mP_{s,x}}\left(G_{s_m}(\omega)H_n(r,\omega_{r})\right)\dif r.
$$
Therefore,
\begin{align*}
\lim_{n\to\infty}
\int^{t_1}_{t_0}\mE^{\mP^n_{s,x}}\left(G_{s_m}(\omega)\mE^{\mP^n_{r,\omega_{r}}}(f(\omega_{t_1}))\right)\dif r
=\int^{t_1}_{t_0}\mE^{\mP_{s,x}}\left(G_{s_m}(\omega)\mE^{\mP_{r,\omega_{r}}}(f(\omega_{t_1}))\right)\dif r,
\end{align*}
which together with \eqref{JG4} and \eqref{JG5} gives \eqref{JG7}.
The proof is complete.
\medskip\\
\begin{proof}[Proof of Theorem \ref{TH11}]
By Lemma \ref{Le34}, we have the existence of $\mP_{s,x}\in\sM^b_{s,x}$, which satisfies the Krylov estimate \eqref{EM9}.
By Lemma \ref{Le36}, we have (i). By Subsection 3.3 we have (ii). 
By Lemma \ref{Le32} and (i), we have (iii).
\end{proof}

\section{Appendix: Properties of space $\widetilde\mH^{\alpha,p}_q$}

In this appendix we prove some important properties about the space $\widetilde\mH^{\alpha,p}_q$.
We need the following lemma, which can be found in \cite[p.205]{Tr} and \cite[Lemma 2.2]{Zh-Zh}.
\bl
\begin{enumerate}[(i)]
\item For any $\alpha\in\mR$ and $p\in(1,\infty)$, there is a $C=C(d,\alpha,p)>0$ such that
\begin{align}\label{EW22}
\|fg\|_{\alpha,p}\leq C\|f\|_{\alpha,p}\|g\|_{|\alpha|+1,\infty}.
\end{align}
\item Let $p\in(1,\infty)$ and $\alpha\in(0,1]$ be fixed. For any $p_1\in[p,\infty)$ 
and $p_2\in[\frac{p_1}{p_1-1},\infty)$ 
with $\frac{1}{p}\leq\frac{1}{p_1}+\frac{1}{p_2}<\frac{1}{p}+\frac{\alpha}{d}$, there is a constant $C>0$ such that
for all $f\in H^{-\alpha,p_1}$ and $g\in H^{\alpha,p_2}$,
\begin{align}\label{EW2}
\|fg\|_{-\alpha,p}\leq C \|f\|_{-\alpha,p_1} \|g\|_{\alpha,p_2}.
\end{align}
\end{enumerate}
\el

The following proposition tells us that the localized norm $\nor\cdot\nor_{\alpha,p;q}$ enjoys the  almost same properties as the global norm $\|\cdot\|_{\alpha,p;q}$.
\bp\label{Pr41}
Let $p,q\in(1,\infty)$ and $\alpha\in\mR$.
\begin{enumerate}[(i)]
\item For $r\not=r'>0$, there is a constant $C=C(d,\alpha,r,r')\geq 1$ such that for all $f\in\widetilde \mH^{\alpha,p}_q$,
\begin{align}\label{GT1}
C^{-1}\sup_{s,z}\|f\chi^{s,z}_{r'}\|_{\alpha,p;q}\leq \sup_{s,z}\|f\chi^{s,z}_r\|_{\alpha,p;q}\leq C \sup_{s,z}\|f\chi^{s,z}_{r'}\|_{\alpha,p;q}.
\end{align}
In other words, the definition of $\widetilde \mH^{\alpha,p}_q$ does not depend on the choice of $r$.
\item Let $(\rho_n)_{n\in\mN}$ be a family of mollifiers in $\mR^d$ and $f_n(t,x):=f(t,\cdot)*\rho_n(x)$. 
For any $f\in\widetilde \mH^{\alpha,p}_q$, it holds that $f_n\in L^q_{loc}(\mR; C^\infty_b(\mR^d))$ and
for some $C=C(d,\alpha,p,q)>0$,
\begin{align}\label{GT2}
\nor f_n\nor_{\alpha,p;q}\leq C\nor f\nor_{\alpha,p;q},\ \forall n\in\mN,
\end{align}
and for any $\varphi\in C^\infty_c(\mR^{d+1})$,
\begin{align}\label{GT3}
\lim_{n\to\infty}\|(f_n-f)\varphi\|_{\alpha,p;q}=0.
\end{align}
\item For any $k\in\mN$, there is a constant $C=C(d,k,\alpha,p,q)\geq 1$ such that for all $f\in \widetilde \mH^{\alpha+k,p}_q$,
$$
C^{-1} \nor f\nor_{\alpha+k,p;q}\leq \nor f\nor_{\alpha,p;q}+\nor \nabla^k f\nor_{\alpha,p;q}\leq C \nor f\nor_{\alpha+k,p;q}.
$$
\item Let $p\in(1,\infty)$ and $\alpha\in(0,1]$, $q\in[1,\infty]$. For any $p_1\in[p,\infty)$ 
and $p_2\in[\frac{p_1}{p_1-1},\infty)$ 
with $\frac{1}{p}\leq\frac{1}{p_1}+\frac{1}{p_2}<\frac{1}{p}+\frac{\alpha}{d}$, and $\frac{1}{q_1}+\frac{1}{q_2}=\frac{1}{q}$, 
there is a constant $C>0$ such that
$$
\nor fg\nor_{-\alpha,p;q}\leq C \nor f\nor_{-\alpha,p_1;q_1} \nor g\nor_{\alpha,p_2;q_2}.
$$
\item $\mL^p_q+\mL^\infty_\infty\subsetneq\widetilde\mL^p_q$.

\end{enumerate}
\ep
\begin{proof}
(i) Let $r>r'$. We first prove the right hand side inequality in \eqref{GT1}. Fix $(s,z)\in\mR^{d+1}$. Notice that
the support of $\chi^{s,z}_r$ is contained in $Q^{s,z}_{2r}$. Clearly, $Q^{s,z}_{2r}$ can be covered 
by finitely many $Q_{r'}^{s_i,z_i}, i=1,\cdots, N$, where $N=N(d,r,r')$ does not depend on $s,z$. 
Let $(\varphi_i)_{i=1}^{N}$ be the partition of unity associated with $\{Q_{r'}^{s_i,z_i}, i=1,\cdots,N\}$ so that
$$
(\varphi_1+\cdots+\varphi_N)|_{Q^{s,z}_{2r}}=1,\ {\rm supp}(\varphi_i)\subset Q_{r'}^{s_i,z_i}.
$$
Thus, due to $\chi^{s_i,z_i}_{r'}|_{Q_{r'}^{s_i,z_i}}=1$, by \eqref{EW22} we have
\begin{align*}
\|f\chi^{s,z}_{r}\|_{\alpha,p;q}&\leq\sum_{i=1}^N\|f\chi^{s,z}_{r}\varphi_i\|_{\alpha,p;q}
=\sum_{i=1}^N\|f\chi^{s_i,z_i}_{r'}\varphi_i\|_{\alpha,p;q}\\
&\leq \sum_{i=1}^N\|f\chi^{s_i,z_i}_{r'}\|_{\alpha,p;q}\|\varphi_i\|_{|\alpha|+1,\infty;\infty}
\leq C\sup_{i=1,\cdots, N}\|f\chi^{s_i,z_i}_{r'}\|_{\alpha,p;q},
\end{align*}
where $C=C(N,\alpha,d,r,r')>0$, which yields the right hand side inequality in \eqref{GT1}.
On the other hand, since $\chi^{s,z}_{r'}=\chi^{s,z}_{2r}\chi^{s,z}_{r'}$, by what we have proved, we have
$$
\|f\chi^{s,z}_{r'}\|_{\alpha,p;q}=\|f\chi^{s,z}_{2r}\chi^{s,z}_{r'}\|_{\alpha,p;q}
\leq \|f\chi^{s,z}_{2r}\|_{\alpha,p;q}\|\chi^{s,z}_{r'}\|_{|\alpha|+1,\infty;\infty}\leq C\|f\chi^{s,z}_{r}\|_{\alpha,p;q},
$$
where $C$ does not depend on $s,z$, which gives the left hand side inequality.

(ii) By the definition of convolutions, it is easy to see that
$$
(\chi^{s,z}_1f_n)(t,x)=\chi^{s,z}_1(t,x)\cdot (f\chi^{s,z}_2)(t,\cdot)*\rho_n(x).
$$
Hence,
$$
\|\chi^{s,z}_1f_n\|_{\alpha,p;q}\lesssim \|\chi^{s,z}_1\|_{|\alpha|+1,\infty;\infty}\|(f\chi^{s,z}_2)_n\|_{\alpha,p;q}
\lesssim \|\chi_1\|_{|\alpha|+1,\infty;\infty}\|f\chi^{s,z}_2\|_{\alpha,p;q},
$$
which gives \eqref{GT2}.
As for \eqref{GT3}, it follows by a finitely covering technique.

(iii) We only prove it for $k=1$. By definition and $\chi^{s,z}_2\nabla\chi^{s,z}_1=\nabla\chi^{s,z}_1$ we have
\begin{align*}
\|(\nabla f) \chi^{s,z}_1\|_{\alpha,p;q}&\leq\|\nabla (f \chi^{s,z}_1)\|_{\alpha,p;q}+\| f \nabla\chi^{s,z}_1\|_{\alpha,p;q}\\
&\lesssim\|f \chi^{s,z}_1\|_{\alpha+1,p;q}+\| f\chi^{s,z}_2\|_{\alpha,p;q}\| \nabla\chi^{s,z}_1\|_{|\alpha|+1,\infty;\infty},
\end{align*}
which in turn gives the right hand side estimate by (i). The left hand side inequality is similar.

(iv) By \eqref{EW2} and $\chi^{s,z}_2\chi^{s,z}_1=\chi^{s,z}_1$, we have
$$
\|(fg)\chi^{s,z}_1\|_{-\alpha,p; q}=\|(f\chi^{s,z}_2) (g\chi^{s,z}_1)\|_{-\alpha,p; q}
\leq\|f\chi^{s,z}_2\|_{-\alpha,p_1;q_1} \|g\chi^{s,z}_1\|_{\alpha,p_2; q_2}.
$$
The desired estimate follows by (i).

(v) Let $\mZ^d$ be the set of all lattice points. Define
$$
f(t,x):=\1_{[0,1]}(t)\sum_{z\in\mZ^d}|x-z|^{-d/p}\1_{|x-z|\leq 1}.
$$
It is easy to see that $f\in \widetilde\mL^p_q$, but $f\notin \mL^p_q+\mL^\infty_\infty$.
\end{proof}

\end{document}